\documentclass[aip,reprint]{revtex4-1}
\usepackage{graphicx}% Include figure files
\usepackage{subfigure}
\usepackage[latin1]{inputenc}
\usepackage{latexsym}
\usepackage{amscd}
\usepackage{amsthm}
\usepackage{amsfonts}
\usepackage{graphicx}
\usepackage{xcolor}
\usepackage{enumerate}
\usepackage{makeidx}
\usepackage{subfigure}
\usepackage{t1enc}
\usepackage{amsmath,amssymb}
\usepackage{ae,aecompl}
\usepackage{booktabs}
\usepackage{verbatim}
\usepackage{lipsum}

\usepackage[abs]{overpic} % to write latex formulas over pictures
\usepackage[font=footnotesize,labelfont=bf]{caption}

\makeatletter
\newtheorem*{rep@theorem}{\rep@title}
\newcommand{\newreptheorem}[2]{%
\newenvironment{rep#1}[1]{%
 \def\rep@title{#2 \ref{##1}}%
 \begin{rep@theorem}}%
 {\end{rep@theorem}}}
\makeatother

\newtheorem{theorem}{Theorem}
\newtheorem{proposition}[theorem]{Proposition}

\newtheorem{lemma}[theorem]{Lemma}
\newtheorem{corollary}[theorem]{Corollary}
\newreptheorem{theorem}{Theorem}

\newtheorem{remark}[theorem]{Remark}
\newtheorem{definition}[theorem]{Definition}

\newcommand{\M}{\text{M}}
\newcommand{\m}{\text{m}}

\newcommand{\NH}{\text{NHIM }}

\usepackage{hyperref}

\begin{document}
The following article has been accepted by Chaos.

After it is published, it will be found at \href{https://publishing.aip.org/resources/librarians/products/journals/}{Link}.

%\preprint{}

\title{Arnold diffusion for an \emph{a priori} unstable Hamiltonian system with 3 $+$ 1/2 degrees of freedom} %Title of paper



\author{A. Delshams}
    %\altaffiliation[Also at ]{Centre de Recerca Matem\`atica (CRM), Barcelona}%Lines break automatically or can be forced with \\
 \email{amadeu.delshams@upc.edu}
 \affiliation{Lab of Geometry and Dynamical Systems and IMTech, Universitat Polit\`ecnica de Catalunya (UPC), 08028 Barcelona, Spain; and Centre de Recerca Matem\`atica (CRM), 08193 Bellaterra, Spain}

\author{A. Granados}%
 \email{albert.granados.corsellas@upc.edu.}
\affiliation{Departament de Matem\`atiques,   Universitat Polit\`ecnica de Catalunya (UPC), 08028 Barcelona, Spain}

\author{R.G. Schaefer}
 \email{rodrigo.goncalves.schaefer@uj.edu.pl.}
\affiliation{Faculty of Mathematics and Computer Science, Jagiellonian University, 30-348 Krak\'ow,Poland}

% Collaboration name, if desird (requires use of superscriptaddress option in \documentclass).

\date{\today}

\begin{abstract}
In the present paper we apply the geometrical mechanism of diffusion in
an \emph{a priori} unstable Hamiltonian system\cite{Chierchia94} with 3 $+$ 1/2 degrees of
freedom.
This mechanism consists of combining iterations of the \emph{inner} and
\emph{outer} dynamics associated to a \emph{Normallly Hyperbolic
Invariant Manifold} (NHIM), to construct diffusing \emph{pseudo-orbits}
and subsequently apply shadowing results to prove the existence of
diffusing orbits of the system.

In addition to proving the existence of diffusion for a wide range of the parameters of the system, an important
part of our study focuses on the search for \emph{Highways}, a particular family of orbits of the outer map (the so-called
\emph{scattering} maps), whose existence is sufficient to ensure a very large drift of the action variables,
with a diffusion time near them that agrees with the optimal
estimates in the literature. Moreover, this optimal diffusion time is calculated, with an explicit calculation of the constants involved.
All these properties are proved by analytical methods and, where necessary, supplemented by numerical calculations.

\end{abstract}


\maketitle %\maketitle must follow title, authors, abstract and \pacs


%==============================================================================
%Introduction
%=============================================================
\begin{quotation}
In this work we study the Arnold diffusion phenomenon on a concrete example of an \emph{a priori} unstable Hamiltonian system with 3+1/2 degrees of freedom.
Thanks to our mechanism, based on iterations of scattering and inner maps, rather than ensuring the existence of Arnold diffusion, we can present examples of diffusing orbits and estimate the diffusing time for some of these orbits.
Throughout this paper, we use analytical methods supplemented by numerical calculations.
\end{quotation}
\section*{Introduction}

In the present paper, we study the geometrical mechanism of diffusion in an \emph{a priori} unstable Hamiltonian system\cite{Chierchia94} with 3 $+$ 1/2 degrees of freedom (d.o.f.)
\begin{equation}\label{eq:ham}
\resizebox{0.9\hsize}{!}{$
H_{\varepsilon}(p,q,I,\varphi,s)=\pm \left(\frac{p^{2}}{2}+\cos q -1\right)+h(I)+\varepsilon f(q) \,g(\varphi,s),
$}
\end{equation}
with $\left(p,q,I,\varphi,s\right)\in \mathbb{R}\times\mathbb{T}\times\mathbb{R}^2\times\mathbb{T}^2\times\mathbb{T}$, where
\begin{align}\label{eq:pert}
h(I)&=h(I_{1},I_{2})=\Omega_{1}I_{1}^{2}/2+
\Omega_{2}I_{2}^{2}/2,\nonumber\\
f(q)& = \cos q,\\% \quad \text{ and }\\
g(\varphi,s)&= g(\varphi_{1},\varphi_{2},s)=a_{1}\cos \varphi_{1} +a_{2}\cos \varphi_{2} +a_{3}\cos s. \nonumber
\end{align}

Combining iterates of the \emph{inner} and the \emph{outer} dynamics associated to the 5D-\emph{Normally Hyperbolic Invariant Manifold} (NHIM)
$\tilde{\Lambda}=\left\{q= 0, p=0\right\}$, to build a diffusing \emph{pseudo-orbit}, and applying shadowing results, we prove the existence of a diffusing orbit of the system.
More precisely, we are able to prove the following theorem on global instability for the following set of parameters:
\begin{equation}\label{eq:parametros}
a_1a_2a_3 \neq 0\quad \text{ and }\quad\left|a_1/a_3\right| + \left|a_2/a_3\right| < 0.625.
\end{equation}
\begin{theorem}\label{thm:main_theorem}
Consider the Hamiltonian \eqref{eq:ham}$+$\eqref{eq:pert}.
Assume $a_1a_2a_3 \neq 0$ and $\left|a_1/a_3\right| + \left|a_2/a_3\right| < 0.625$.
Then, for every $0<\delta < 1$ and $R >0$,  there exists $\varepsilon_0=\varepsilon_0(\delta,R)>0$ such that for every $I_+, I_-$ satisfying $\left|I_{\pm}\right|< R$,
there exists an orbit $\tilde{x}(t)$ and $T >0$, such that
\begin{align*}
\left|I(0) - I_- \right| \leq \delta\quad\text{and} &\quad \left|I(T) - I_+ \right| \leq \delta.
\end{align*}
\end{theorem}

Arnold diffusion\cite{Arnold64}, or global instability in nearly-integrable Hamiltonian systems, has been revitalized in the last decade, thanks to the appearance of very important results for a priori stable general Hamiltonians\cite{BernardKZ16,KaloshinK20,GideaM22,ChengX23}. This work is part of the study of a priori unstable Hamiltonians, in which, thanks to the use of geometrical methods\cite{Delshams2000,Seara2006,GideaLS20}, it is possible to design concrete paths of unstable trajectories, and even measure the time spent along such unstable trajectories.

It is worth remarking that the notion of a priori unstable system was introduced in the seminal paper~\cite{Chierchia94}, along with geometric methods following the ideas of Arnold's article~\cite{Arnold64}. The geometrical method~\cite{Seara2006} used in this paper is based on combining two different dynamics on a normally hyperbolic manifold (NHIM) of the system, the inner map and the outer or scattering maps, to produce unstable pseudo-trajectories (or diffusive trajectories, as they are also called). The design of diffusive trajectories essentially depends on the scattering maps, although to apply it properly, the inner map must be taken into account.

This paper is dedicated to showing how, considering perhaps the simplest perturbation based solely on three parameters corresponding to three harmonics, one can prove the existence of Arnold diffusion in a very large open of the parameter set. And in addition, fast diffusion trajectories (which we call \emph{Highways}) can be designed, with a quantitative estimate of the time used. In this sense, it should be noted that we are not talking about general or generic perturbations, but about concrete perturbations where we apply our results. Indeed, a proof that this geometric method is useful for applications is the global instability results in Celestial Mechanics~\cite{Delshams16,Fejoz16,Capinski2017,Clarke22} that have been obtained recently.

The case of an \emph{a priori} unstable Hamiltonian system with 2$+$1/2 degrees of freedom and for a `complete' family of perturbations, that is, based on two arbitrary independent harmonics, was already dealt by the authors\cite{Delshams2017,Delshams2018}.
Although there are some similarities between the mechanism for 3$+$1/2 and 2$+$1/2 degrees of freedom, there are also very remarkable differences.

One of these differences appears in our main tool, the \emph{scattering maps} defined on the NHIM $\tilde{\Lambda}$.
As we explain in Section~\ref{sec:scattering_map}, the trajectories of a scattering map $\mathcal{S}$ are given, up to order $\mathcal{O}(\varepsilon^2)$, by the $-\varepsilon$-time flow of a Hamiltonian $\mathcal{L}^*(I,\theta = \varphi - Is)$, called the \emph{reduced Poincar\'e function}.
In 2$+$1/2 d.o.f., we took advantage of the fact that the Hamiltonian system associated to $\mathcal{L}^*(I, \theta)$ had 1 degree of freedom, so it was integrable \cite{Delshams2017,Delshams2018}.
This integrability allowed us to describe the scattering maps and their bifurcations entirely by looking at the behavior of the level curves of $\mathcal{L}^*$.
In contrast, for the Hamiltonian (1)$+$(2) considered in this work, $\mathcal{L}^*(I,\theta)$ has 2 d.o.f. and, in general, is not integrable (see Fig.~5 for an illustration).
Consequently, we will not perform a complete description of all diffusion paths connecting two arbitrary values of the actions, but will only look for some suitable diffusion paths.
It is very important to emphasize that it is not necessary to know \textbf{all} the outer dynamics (the scattering map) generated by $\mathcal{L}^*(I,\theta)$ to produce global instability. It is enough to locate the trajectories of the scattering map that produce the most amount of diffusion in the action variables $I$. This is one of the advantages of having separated the dynamics near the NHIM into two different dynamics.
Therefore, to show fast diffusion paths, we will focus a significant part of our study near \emph{Highways}.

Highways were introduced\cite{Delshams2017}, in the case of 2$+$1/2 d.o.f., as two curves contained in a specific level $\mathcal{L}^*= C_{\text{h}}$ of the reduced Poincar\'e function which had excellent properties to perform fast diffusion: the Highways were curves $\theta = \Theta(I)$, so they were ``vertical'' on the plane $(\theta = \varphi - Is , I)$, which means that the action $I$ increases or decreases significantly along the Highways.
The global or local existence of Highways depended on the perturbation parameters\cite{Delshams2017}.
Thanks to the values of the parameter~\eqref{eq:parametros} that we consider, we will show in Theorem \ref{thm:global_existence_highways} that the Highways are globally defined for systems \eqref{eq:ham}$+$\eqref{eq:pert}.
In Section \ref{sec:global_instability}, we provide a more geometrical definition for Highways: unlike the 2 level curves in the case of 2$+$1/2 d.o.f., there are now four
Highways $\theta=\Theta(I)$ which are \emph{Lagrangian} surfaces contained in a specific 3D manifold $\mathcal{L^*}= C_{\text{h}}$.
In addition, we will provide an estimate $T_{\mbox{d}}=T_{\mbox{s}}/\varepsilon\left(\log(C/\varepsilon)+\mathcal{O}(\varepsilon^{b})\right)$, $0<b<1$, of the time of diffusion for orbits close to Highways, which coincides with previous general optimal estimates~\cite{Berti2003,Cresson2003,Tre04}, but includes an explicit expression of $T_s$ \eqref{eq:formula_time} as well as of $C$.
As will be verified, this estimate is similar to the one found in Ref.~\onlinecite{Delshams2017} and consists mainly of the time under the scattering map.
The dynamics on these Highways will be described with the help of numerical computations.
In passing, it is worth noticing that there are no Highways in general\cite{Delshams2018}.

The inner dynamics in $\tilde{\Lambda}$ is described in Section~\ref{sec:inner}.
In general, a Hamiltonian system with 3$+$1/2 d.o.f. restricted to a 5D $\tilde{\Lambda}$ gives rise to a Hamiltonian $K_{\varepsilon}=H_{\varepsilon}|_{\tilde{\Lambda}}$ with 2$+$1/2 degrees of freedom where double resonances can appear.
Due to the hypotheses in \eqref{eq:ham}$+$\eqref{eq:pert}, there exists only one double resonance at $ I = 0$ and the Hamiltonian $K_{\varepsilon}$ that governs the inner dynamics is integrable.
We observe that the inner dynamics in Ref.~\onlinecite{Delshams2017,Delshams2018} have different features.
While in Ref.~\onlinecite{Delshams2017} it is integrable and has only one resonance of first order at $I = 0$, in Ref.~\onlinecite{Delshams2018} the inner dynamics is no longer integrable and has two resonances of first order, $I= 0$ and $I= 1$.

The system given by Hamiltonian \eqref{eq:ham}$+$\eqref{eq:pert} is a direct generalization of the Hamiltonian considered in Ref.~\onlinecite{Delshams2017}.
It is worth noting that in the Hamiltonian systems with 2$+$1/2 degrees of freedom of Ref.~\onlinecite{Delshams2017,Delshams2018}, there existed remarkable curves called \emph{crests} that played an essential role in understanding the bifurcation of scattering maps.
The crests (or ridges) were introduced in Ref.~\onlinecite{Delshams2011} and studied in Ref.~\onlinecite{Delshams2017,Delshams2018}.
In Ref.~\onlinecite{Delshams2017}, the crest $\mathcal{C}(I)$ is formed simply by two curves, $\mathcal{C}_{\M}(I)$ and $\mathcal{C}_{\m}(I)$, in the plane $\left(\varphi, s\right)$.
On each curve, there exists at least one locally defined scattering map.
The shape of these curves depends on the value of the coefficients $a_i$ of the perturbation and has a direct influence on the behavior and domain of the scattering maps.
In Section~\ref{sec:scattering_map}, we study the bifurcation of the crests of Hamiltonian \eqref{eq:ham}$+$\eqref{eq:pert} with respect to parameters $a_1$, $a_2$ and $a_3$.
In this case, the crest $\mathcal{C}(I)$ is formed either by \emph{two} surfaces or by only \emph{one} surface in the plane $\left(\varphi_1, \varphi_2, s\right)$ .
In this paper, we restrict ourselves to the case where there are only two scattering maps associated to the crests and are globally defined, that is, we have the scattering maps $\mathcal{S}_{\M}$ and $\mathcal{S}_{\m}$ associated to $\mathcal{C}_{\M}(I)$ and $\mathcal{C}_{\m}(I)$, respectively.
We show that this happens just for condition \eqref{eq:parametros}.
This emphasizes the similarity of \eqref{eq:ham}$+$\eqref{eq:pert} with the system in Ref.~\onlinecite{Delshams2017}.

Our model describes the dynamics of one pendulum plus two rotors, but we could also consider several pendula plus two rotors.
It is worth remarking that in a system with more pendula, there are generally more homoclinic manifolds with respect to its NHIM.
Therefore, more distinct scattering maps can be defined, providing more possibilities to detect global instability.

In Ref.~\onlinecite{Delshams2018}, it was observed that choosing a perturbation function of the form $f(q)g(\varphi,s)$ simplifies the calculation of the
 Poincar\'e-Melnikov potential~\eqref{eq:mel_pot_gen}. However, it is worth noting that it is possible to readily extend Theorem~\ref{thm:main_theorem} to the case where $f(q)g(\varphi,s)$ is a trigonometric or meromorphic polynomial in $q$, although a more detailed analysis may be necessary for the Hamiltonian of the inner dynamics.

This work is related to Ref.~\onlinecite{DelshamsLS07b}, but the philosophy is quite different.
The approach in Ref.~\onlinecite{DelshamsLS07b} is more general and theoretical, and the hypotheses depended on several conditions that had to be tested.
In this paper, we can prove the existence of global instability for a `\emph{complete}' family of perturbations \eqref{eq:pert}$+$\eqref{eq:parametros} depending on three harmonics.
We are also interested in computing \emph{explicitly} some fast diffusing orbits.
We apply as much as possible analytical tools and, when necessary, complete our analytical results with numerical computations.

The system discussed in the present paper is a particular case of Hamiltonian \eqref{eq:ham} with a function $g(\varphi,s)$ satisfying
\begin{equation*}
g(\varphi,s)=a_{1}\cos \varphi_{1} +a_{2}\cos \varphi_{2} +a_{3}\cos(k\cdot\varphi-s),
\end{equation*}
with $k = (k_1 , k_2) \in \mathbb{Z}^2$.
Indeed, the perturbation considered here is just $k = 0$, which is a direct generalization of the system studied in Ref.~\onlinecite{Delshams2017}.
The next step that we want to address in a future work is the case $k\neq 0$, which would be a generalization of the system studied in Ref.~\onlinecite{Delshams2018}, and in particular the case $k = (1,1 )$, which was used in Ref.~\onlinecite{DelshamsLS07b} to illustrate the results obtained there.
To deal with this case as well would make this paper longer and somewhat more technical. Although our method can also be applied, understanding the regions where the scattering maps are well defined is more complicated due to the more complex behavior of the crests and requires more detailed study.

There are several works in the literature dealing with Arnold diffusion in similar or more general a priori unstable models. For example, in Ref.~\onlinecite{Chen2022}, the authors proved the existence of Arnold diffusion in a much more general context. In Ref.~\onlinecite{Gallavotti2000}, a model with different time scales is considered. A very similar perturbation function for 2$+$1/2 degrees of freedom is studied in Ref.~\onlinecite{Akingbade23}, but for an unperturbed dissipative part, and a Hamiltonian with two pendula and one rotor is considered in \onlinecite{Bessi1997}.
Different geometrical and variational methods appear also in several other papers~\cite{Berti2003,Tre04,GelfreichSV13,BernardKZ16,Davletshin2016}.
The novelty of our work is that we not only prove the existence of diffusion through abstract reasoning, but also show the diffusion paths for concrete Hamiltonian systems. Furthermore, we present a method to identify \emph{fast} diffusion paths.

%================================================================================================
% UNPERTURBED CASE
%=====================================================================================================
\section{Unperturbed case}
%\label{sec:unperturbed_case}

For $\varepsilon = 0$, Hamiltonian \eqref{eq:ham}$+$\eqref{eq:pert} becomes
$$H_0(p,q, I ,\varphi,s) = \pm \left(\frac{p^{2}}{2}+\cos q -1\right)+\frac{\Omega_1 I_1^2 }{2} + \frac{\Omega_2 I_2^2}{2},$$
with associated equations
\begin{align*}
&\dot{q} = p& &\dot{p} = \sin q \\
&\dot{\varphi_{1}} = \Omega_1 I_1=: \omega_{1}&  &\dot{I_{1}}= 0  \nonumber\\
&\dot{\varphi_{2}} = \Omega_2I_2=: \omega_{2}&  &\dot{I_{2}}= 0  \nonumber \\
&\dot{s}=1 , \nonumber
\end{align*}
This system consists of one pendulum plus a $2$ d.o.f rotor.
From the equations above, $I_{1}$ and $I_{2}$ are constants of motion, and the flow has the form
$$\Phi_{t}(p,q, I ,\varphi) =(p(t),\,q(t),I ,\varphi + t\omega),$$
where $\omega = (\omega_1 , \omega_2):= (\Omega_1I_1 , \Omega_2 I_2)$.
To include the frequency of the time, we also will use the frequency vector $\tilde{\omega} = (\omega_1, \omega_2, 1)$.

Observe that $(p_0,q_0)=(0,0)$ is a saddle point on the plane $(p,q)$ with unstable and stable invariant curves.
These invariant curves coincide and separate the behavior of orbits and, for this reason, are called \emph{separatrices}.
In addition, they can be parametrized by
\begin{equation}
\label{eq:separatrix}
(p_0(t), q_0(t)) = \left(\frac{\pm 2}{\cosh t} , 4\arctan e^{\pm t}\right).
\end{equation}

For any $I\in\mathbb{R}^2$,
$ \mathcal{T}_{I} = \{(0,0,I,\varphi,s); \varphi,s \in \mathbb{T}^{3}\}$
is an invariant torus under the flow of the system with the frequency $\omega = \left(\Omega_1I_1 , \Omega_2 I_2\right)$ and is called a \emph{whiskered torus}.
For each whiskered torus, we have associated coincident stable and unstable manifolds called \emph{whiskers}, which we denote by
\begin{equation*}
W^0\mathcal{T} = \left\{(p_0(\tau) , q_0(\tau) , I, \varphi, s):\tau\in\mathbb{R}, (\varphi , s)\in\mathbb{T}^2\right\}.
\end{equation*}

The union of all whiskered tori $\mathcal{T}$
\begin{equation}\label{eq:NHIM_man}
\tilde{\Lambda} = \{(0,0,I,\varphi,s):(I,\varphi,s)\,\in\,\mathbb{R}\times\,\mathbb{T}^{3}\}
\end{equation}
is a 5D-\emph{Normally Hyperbolic Invariant Manifold} (NHIM) with 6D-coincident stable and unstable invariant manifolds given by
\begin{equation*}
W^0\tilde{\Lambda} = \{(p_0(\tau),q_0(\tau),I,\varphi,s):\tau\in\mathbb{R}, (I,\varphi,s)\,\in\,\mathbb{R}^2\times\,\mathbb{T}^{3}\}.
\end{equation*}

For $\varepsilon\neq 0$, there exists a NHIM $\tilde{\Lambda}_{\varepsilon}$, which in our system, due to the fact that $f(q) =\cos q$ in \eqref{eq:ham}, coincides with the NHIM $\tilde{\Lambda}$.
Nevertheless, its stable manifold $W^{\text{s}}(\tilde{\Lambda}_{\varepsilon})$ and unstable manifold $W^{\text{u}}(\tilde{\Lambda}_{\varepsilon})$ no longer coincide, that is, the separatrices split.
%===========================================================
% INNER DYNAMICS
%===========================================================
\section{Inner dynamics}
\label{sec:inner}

The inner dynamics is derived from the restriction of the Hamiltonian~\eqref{eq:ham} to $\tilde{\Lambda}$ given in \eqref{eq:NHIM_man}, that is,

\begin{align*}
K_{\varepsilon}( I , \varphi , s ) = &\sum\limits_{k=1}^2\frac{\Omega_kI_k^2}{2} + \varepsilon\left(\sum\limits_{k=1}^2 a_k\cos\varphi_k + a_3\cos s\right)
\end{align*}%\label{eq:ham_inner}
and its associated equations are
\begin{align}\label{eq:inner_eq}
\dot{\varphi}_1 &= \omega_1 = \Omega_1 I_1 & \dot{I}_1 &= \varepsilon a_1\sin\varphi_1 \\
\dot{\varphi}_2 &= \omega_2 = \Omega_2 I_2 & \dot{I}_2 &= \varepsilon a_2\sin\varphi_2 \nonumber\\
\dot{s}& = 1. & & \nonumber
\end{align}

\begin{figure}[h]
\centering
\includegraphics[scale=0.22]{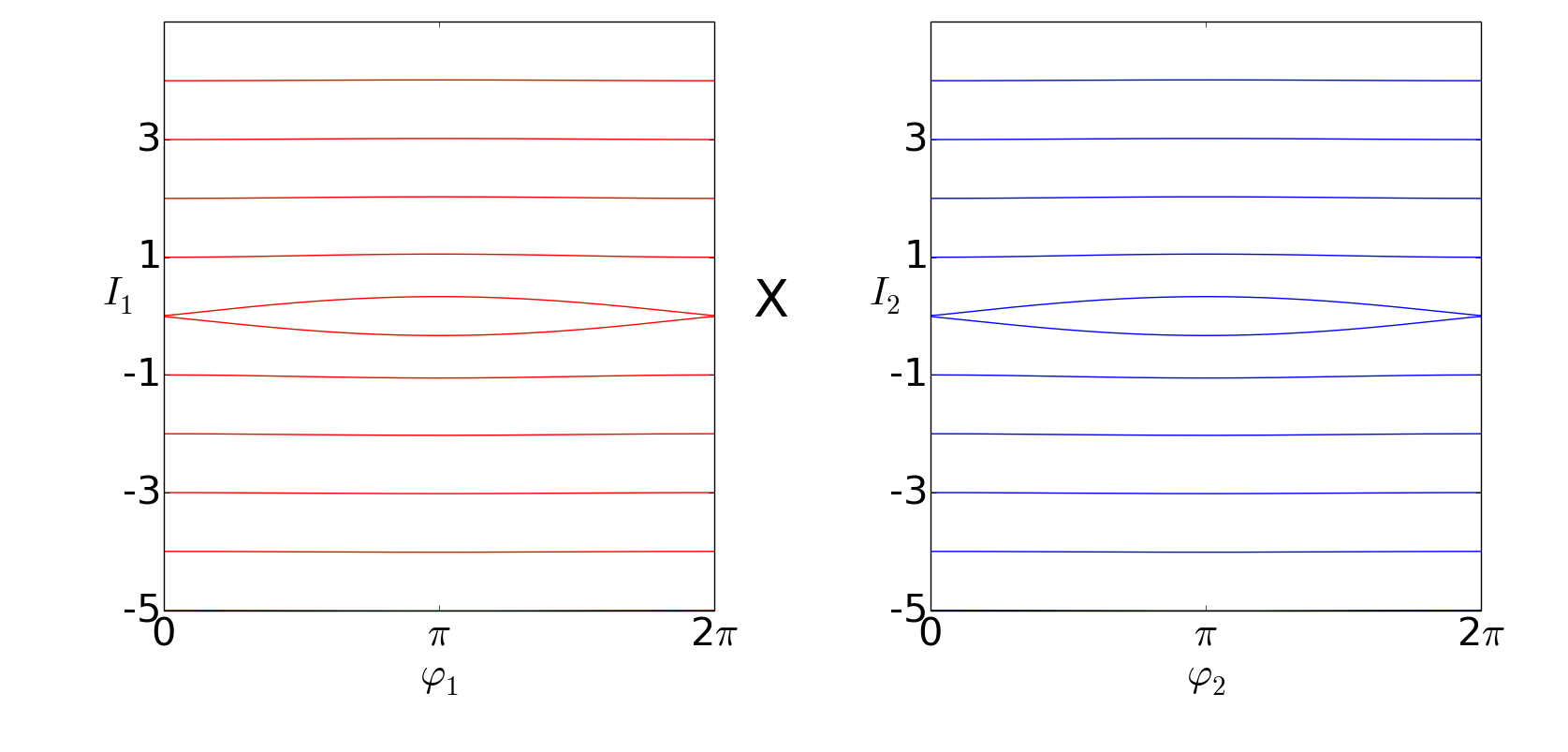}
\caption{Inner dynamics}
\label{fig:inner_map}
\end{figure}

Note that the inner dynamics is integrable, with first integrals
\begin{align*}
F_1(I_1,\varphi_1) &= \frac{\Omega_1I_1^2}{2} + \varepsilon a_1\left(\cos\varphi_1 - 1\right),\\
F_2(I_2,\varphi_2) &= \frac{\Omega_2I_2^2}{2} +\varepsilon a_2\left(\cos\varphi_2 - 1\right),
\end{align*}
in involution.
The inner dynamics is described in Fig.~\ref{fig:inner_map}, so there are two resonances, one centered at $I_1 = 0$ and other at $I_2 = 0$.

%=============================== figure: horizontal crest
%=========================================================
\begin{remark}\label{rem:inner_horizontal}
When $\varepsilon$ is small enough, the level curves of $F_1$ and $F_2$ are almost \emph{flat} or \emph{horizontal} in the action $I = (I_1,I_2)$, i.e., the values of $I_1$ and $I_2$ remain almost constant, see Fig.~\ref{fig:inner_map}.
Since there is a continuous foliation of invariant tori, a genuine `large gap problem'~\cite{Seara2006} does not appear.	
\end{remark}

\begin{remark}
There is a double resonance at $I_1 = I_2 = 0$.
The study of dynamics close to double resonances is a challenging problem and is out of the scope of this work. Since, in our case, the double resonance is just the point $I= (0,0)$, we will avoid it when necessary.
\end{remark}

%==========================================================================================================
% SCATTERING MAP
%======================================================================================================
\section{Scattering map}
\label{sec:scattering_map}
%----------------------------------------------------------------------------------------------------
% definition of scattering map
%----------------------------------------------------------------------------------------------------
\subsection{Definition of scattering map }
%\label{subsec:def_of_scattering_map}

The notion of a scattering map on a \NH was introduced by Delshams et al.\cite{Delshams2000}. It plays a central r\^ole in our mechanism for detecting diffusion.
Let $W$ be an open set of $\left[-I_1^* , I_1^*\right]\times \left[-I_2^* , I_2^*\right]\times \mathbb{T}^3$ such that the invariant manifolds of \NH $\tilde{\Lambda}$
introduced in~\eqref{eq:NHIM_man} intersect transversely along a homoclinic manifold $\Gamma~=~\left\{\tilde{z}(I,\varphi,s;\varepsilon) , (I,\varphi,s)\in W\right\}$ and for any $\tilde{z}\in\Gamma$ there exists a unique $\tilde{x}_{\pm} = \tilde{x}_{\pm}(I,\varphi,s;\varepsilon)\in\tilde{\Lambda}$ such that $\tilde{z}~\in~W_{\varepsilon}^{s}(x_-)\cap W_{\varepsilon}^{u}(\tilde{x}_+)$.
Let
$$H_{\pm} = \bigcup\left\{\tilde{x}_{\pm}(I ,\varphi , s ; \varepsilon) : (I,\varphi,s)\in W\right\}.$$
The scattering map associated to $\Gamma$ is the map
\begin{eqnarray*}
S: H_- & \longrightarrow & H_+\\
\tilde{x}_- &\longmapsto & S(\tilde{x}_-) = \tilde{x}_+.
\end{eqnarray*}

Notice that the domain of definition of the scattering map depends on the homoclinic manifold chosen.
Therefore, for the characterization of the scattering maps, it is required to select the homoclinic manifold $\Gamma$, which can be done using the Poincar\'{e}-Melnikov theory. We have the following proposition\cite{Seara2006,Delshams2011}

\begin{proposition}
Given $(I,\varphi,s)\,\in\,\left[-I_1^{*},I_1^{*}\right]\times \left[-I_1^{*},I_1^{*}\right]\,\times\,\mathbb{T}^{3}$, assume that the real function
\begin{equation}\label{eq: SM_critical_point}
\tau\,\in\,\mathbb{R}\,\longmapsto\,\mathcal{L}(I ,\varphi -\tau\omega ,\,s-\tau)\,\in\,\mathbb{R}
\end{equation}
has a non-degenerate critical point $\tau^{*}\, =\, \tau^*(I,\varphi,s)$, where $\omega = (\omega_1 , \omega_2)= \left(\Omega_1I_1 , \Omega_2I_2\right)$ and
\begin{equation*}
\mathcal{L}(I , \varphi,s):=\int_{-\infty}^{+\infty}\left(f(q_{0}(\rho)) - f(0)\right)g(\varphi+ \rho\omega,s+\rho;0)d\rho.
\end{equation*}
Then, for $0\,<\,\varepsilon$ small enough, there exists a unique transverse homoclinic point $\tilde{z}$ to $\tilde{\Lambda}_{\varepsilon}$ of Hamiltonian~\eqref{eq:ham}, which is $\varepsilon$-close to the point
$\tilde{z}^{*}(I,\varphi,s)\,=\,(p_{0}(\tau^{*}),q_{0}(\tau^{*}),I,\varphi,s)\,\in\,W^{0}(\tilde{\Lambda})$:
%\label{eq:heteroclinic_point}
\begin{align*}
&\tilde{z}=\tilde{z}(I,\varphi,s)=(p_{0}(\tau^{*})+O(\varepsilon), q_{0}(\tau^{*})+O(\varepsilon),I,\varphi,s)\\
&\text{and }\tilde{z}\in\,W^{u}(\tilde{\Lambda_ {\varepsilon}})\,\pitchfork\,W^{s}(\tilde{\Lambda_{\varepsilon}}).
\end{align*}
\end{proposition}

The function $\mathcal{L}$ is called the \emph{Melnikov potential} of Hamiltonian \eqref{eq:ham}, and using~\eqref{eq:pert} and \eqref{eq:separatrix} takes the form
\begin{equation}\label{eq:mel_pot_gen}
\mathcal{L}(I,\varphi,s) = A_1\cos\varphi_1  + A_2\cos\varphi_2 + A_3\cos s,
\end{equation}
where
\begin{eqnarray}
A_i:=A(\tilde{\omega}_i,a_i)= \frac{2\pi\tilde{\omega}_i a_i}{\sinh(\pi\tilde{\omega}_i/2)},\quad i = 1,2,3
\label{eq:coef_A}
\end{eqnarray}
and $\tilde{\omega} = (\omega_1,\omega_2,1 )$, for $\omega_i \neq 0$, and $A_i = 4a_i$ for $\omega_i= 0$.
The homoclinic manifold $\Gamma$ is determined by the function $\tau^*(I ,\varphi ,s)$.
Once a function $\tau^*(I,\varphi,s)$ is chosen, by the geometric properties of the scattering map, see \cite{Delshams2008,Delshams2009,Delshams2011}, the scattering map $S = S_{\tau^*}$ has the explicit form
\begin{equation*}
S(I,\varphi,s) = \left(I \ + \varepsilon\nabla_{\varphi}L^* + \mathcal{O}(\varepsilon^2), \varphi - \varepsilon\nabla_{I} L^* +\mathcal{O}(\varepsilon^2), s \right),
\end{equation*}
where
\begin{equation}\label{eq:L^*-def}
\resizebox{0.9\hsize}{!}{$
L^* = L^*(I,\varphi,s) = \mathcal{L}\left(I,\varphi -\tau^*(I,\varphi,s)\omega , s-\tau^*(I,\varphi,s)\right).$}
\end{equation}

Note that the variable $s$ is fixed under the scattering map.
As a consequence, introducing the variable
\begin{equation}
\theta = \varphi - s\,\omega,\label{eq:def_theta}
\end{equation}
we can define the \emph{reduced Poincar\'{e} function} as
\begin{equation}\label{eq:reduced_poincare_function}
\mathcal{L}^{*}(I,\theta) := L^*(I,\varphi - s\omega  , 0) = L^*(I,\varphi,s).
\end{equation}
In these variables $(I,\theta)$ the scattering map has the simple form
\begin{align}
\mathcal{S}(I,\theta) =& \left( I + \varepsilon\frac{\partial \mathcal{L}^*}{\partial\theta}(I,\theta) + \mathcal{O}(\varepsilon^2) ,\right.\nonumber\\
&\left. \theta - \varepsilon\frac{\partial \mathcal{L}^*}{\partial I}(I,\theta) + \mathcal{O}(\varepsilon^2)   \right).\label{eq:form_scatteringmap}
\end{align}
So, up to $\mathcal{O}(\varepsilon^2)$ terms, $\mathcal{S}(I,\theta)$ is the $-\varepsilon$ times flow of the \emph{autonomous} Hamiltonian $\mathcal{L}^*(I,\theta)$.
In particular, a finite number of iterates under the scattering map follow the level curves of $\mathcal{L}^*$ up to $\mathcal{O}(\varepsilon^2)$.

It is easy to see\cite{Delshams2011,Delshams2018} that the reduced Poincar\'e function $\mathcal{L}^*(I,\theta)$ in \eqref{eq:reduced_poincare_function} is equivalent to
$$\mathcal{L}^*(I,\theta) ~=~ \mathcal{L}(I,\theta -~ \tau^*(I,\theta)\omega, -\tau^*(I,\theta)),$$
 where
\begin{align}\label{eq:tau_theta}
\tau^*(I,\theta) = \tau^*(I,\varphi,s) - s.
\end{align}
Therefore, from \eqref{eq:mel_pot_gen} the reduced Poincar\'e function $\mathcal{L}^*(I,\theta)$ can be explicitly expressed as
\begin{align}\label{eq:reduced_poincare_explicit}
\mathcal{L}^*(I,\theta) &= \sum\limits_{k=1}^2A_k\cos(\theta_k - \omega_k\tau^*(I,\theta))\nonumber\\
&+ A_3\cos(-\tau^*(I,\theta)).
\end{align}
Along this paper, both $\tau^*(I,\varphi,s)$ and $\tau^*(I,\theta)$ will be used at our convenience.
It is important to note that, as the variable $s$ is fixed for $\mathcal{S}(I,\theta)$, it plays the r\^ole of a parameter, which for simplicity\cite{Delshams2017} will be taken $s = 0$.

%--------------------------------------------------------------------------------
% Crests and {\NH} lines
%-----------------------------------------------------------------------------
\subsection{Crests and {\NH} lines}

We have seen that the function $\tau^*$ plays a key role in our study.
Therefore, we are interested in finding the critical points $\tau^* = \tau^*(I,\varphi,s) $ of function \eqref{eq: SM_critical_point} or, for our concrete case \eqref{eq:mel_pot_gen}, $\tau^*$ solution of
\begin{equation}\label{eq:tau^*}
\resizebox{0.88\hsize}{!}{$
\frac{\partial \mathcal{L}}{\partial \tau}(I,\varphi-\omega \tau , s-\tau) = \nabla_{(\varphi,s)}\mathcal{L}(I,\varphi-\omega\tau,s-\tau)\cdot\tilde{\omega}= 0.$}
\end{equation}
This equation can be seen from two equivalent geometrical viewpoints.
The first one is that to find $\tau^* = \tau^*(I,\varphi,s)$ satisfying \eqref{eq:tau^*} for any $(I,\varphi,s)\in \left[-I_1^* , I_1^* \right]\times\left[-I_2^* , I_2^*\right]\times\mathbb{T}^3$ is the same as looking for the extrema of $\mathcal{L}$ on the \emph{{\NH} line}
\begin{equation}
R(I,\varphi,s) = \left\{(I,\varphi - \tau\omega , s-\tau) : \tau \in \mathbb{R}\right\}.\label{eq:nhim_line}
\end{equation}
The other point of view is that fixing $(I,\varphi,s)$, a solution $\tau^*$ of \eqref{eq:tau^*} is equivalent to finding intersections between a \NH line \eqref{eq:nhim_line} and a set of points determined by the equation $\nabla_{(\varphi,s)}\mathcal{L}(I,\varphi,s)\cdot\tilde{\omega} =0.$

\begin{remark}
Note that $R(I,\varphi,s)$ on $\left\{I\right\}\times \mathbb{T}^3$ is either a closed line when $\omega\in\mathbb{Q}^2$, or densely fills a torus of dimension $2$ or $3$.
\end{remark}

\begin{remark}
By taking $\theta$ as defined in \eqref{eq:def_theta} and $s = 0$, we can rewrite $R(I,\varphi,s)$ on variables $(I,\theta)$, that is, $R(I,\varphi,s)$ is equivalent to
\begin{equation}
R(I,\theta) = \left\{(I,\theta - \tau\omega , -\tau) : \tau \in \mathbb{R}\right\}.
\end{equation}

\end{remark}

\begin{definition}\cite{Delshams2011}\label{def:crest}
 A \emph{crest} or \emph{ridge} $\mathcal{C}(I)$ is the set of points $\left\{(I,\varphi,s),\, (\varphi,s)\in \mathbb{T}^3\right\}$ such that
\begin{equation*}
\frac{\partial \mathcal{L}}{\partial \tau}(I , \varphi - \tau\omega ,s - \tau)|_{\tau = 0} = 0,
\end{equation*}
or equivalently,
\begin{equation}
\nabla_{(\varphi,s)}\mathcal{L}(I,\varphi,s)\cdot\tilde{\omega} =0.\label{eq:def_crests}
\end{equation}
\end{definition}
Fixing $I$, Eq.~\eqref{eq:def_crests} defines, at least locally, a surface in the variables $(\varphi_1,\varphi_2, s)\in \mathbb{T}^3$, that we want to characterize.

From the expression \eqref{eq:mel_pot_gen} of the Melnikov Potential $\mathcal{L}(I,\varphi,s)$, Eq.~\eqref{eq:def_crests} can be rewritten as
\begin{equation}
\alpha(\omega_1)\mu_1 \sin\varphi_1 + \alpha(\omega_2)\mu_2\sin\varphi_2 + \sin s = 0\label{eq:crests_cms}
\end{equation}
where, for $i = 1,2$,
\begin{equation} \label{eq:def_alpha}
\mu_i =\frac{a_i}{a_3} \quad \text{ and }\quad \alpha(\omega_i) = \left(\omega_i\right)^{2}\frac{\sinh(\pi/2)}{\sinh(\omega_i\pi/2)}.
\end{equation}
Observe that $\alpha$ is well defined for any value of $\omega_i$.
For any fixed $I$, to understand the intersection between \NH lines $R(I,\cdot)$ and the crest $\mathcal{C}(I)$, we first need to study how the surfaces contained in the crest look like for different values of $\mu_i=a_i/a_3$ and $\omega_i$, for $i = 1,2$.

\begin{remark}
We wish to emphasize the similarity between Eq.~\eqref{eq:crests_cms} of the crest with the equation of the crests studied for 2$+$1/2 d.o.f.~\cite{Delshams2017,Delshams2018}
\end{remark}

Eq.~\eqref{eq:crests_cms} is just a linear equation for the variables $\sin\varphi_1$, $\sin\varphi_2$ and $\sin s$.
To parametrize the crest, we want to isolate one of these variables with respect to the other two.
We begin with the case
\begin{equation}\label{eq:horizontal_cond}
\text{a)}\quad\left|\alpha(\omega_1)\mu_1\sin\varphi_1 + \alpha(\omega_2)\mu_2\sin\varphi_2\right| \leq 1,
\end{equation}
where we can write $s$ as a function of $\varphi_1$ and $\varphi_2$ for any $\left(\varphi_1 , \varphi_2\right)\in\mathbb{T}^2$, more precisely, we have the two functions
\begin{align*}
s = \left\{\begin{matrix}
&\xi_{\M}(I,\varphi) :=\arcsin\left( \alpha(\omega_1)\mu_1\sin\varphi_1\right.\\
&\left. + \alpha(\omega_2)\mu_2\sin\varphi_2\right) \mod 2\pi \vspace{0.3cm}\\
&\xi_{\m}(I,\varphi) :=-\arcsin\left( \alpha(\omega_1)\mu_1\sin\varphi_1\right.\\
&\left. + \alpha(\omega_2)\mu_2\sin\varphi_2\right) + \pi \mod 2\pi .
\end{matrix}\right.
\end{align*}

\begin{figure}[h]
\centering
\includegraphics[scale=0.4]{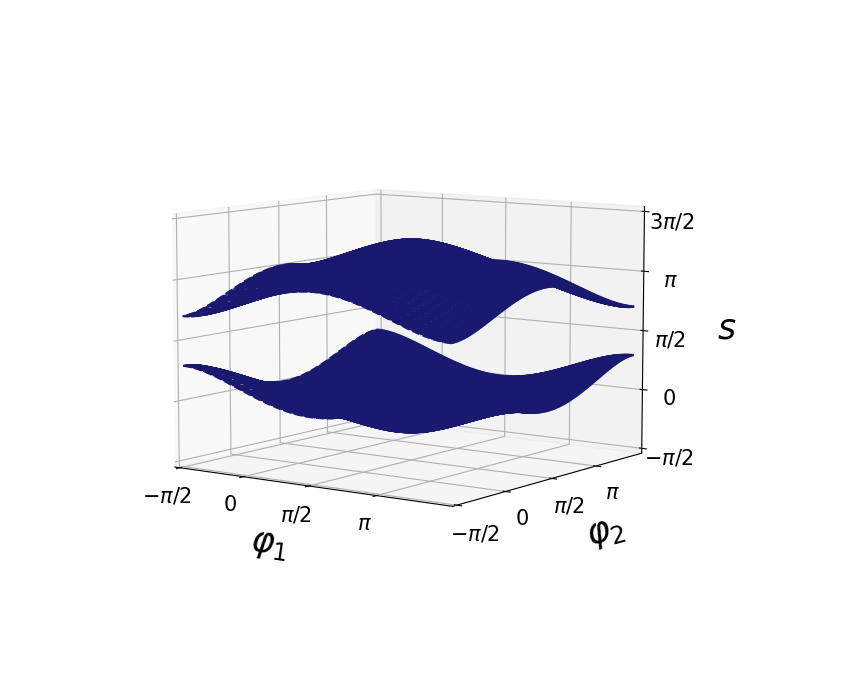}
\caption{\emph{Horizontal Crest} $\mathcal{C}(I)$ formed by two surfaces: $\mu_1 = \mu_2 = 0.4$ and $\omega_1 = \omega_2 = 1$. We plotted $\varphi_1,\varphi_2,s\in \left[-\pi/2, 3\pi/2\right)$ for a better illustration.}
\label{fig:cristas_horizontais}
\end{figure}

Then the crest $\mathcal{C}(I)$ is formed by two surfaces $\mathcal{C}_{\text{M}}(I)$, $\mathcal{C}_{\text{m}}(I)$ which are simply the \emph{global graphics} in the space $\left(\varphi_1, \varphi_2, s\right)$ of the functions $\xi_{\text{M}}(I,\varphi)$, $\xi_{\text{m}}(I,\varphi)$, defined for all $\varphi\in\mathbb{T}^2$.
According to the notation used for 2$+$1/2 d.o.f.\cite{Delshams2017,Delshams2018}, the crest $\mathcal{C}(I)$ is called a \emph{horizontal} crest; see Fig.~\ref{fig:cristas_horizontais}.
From the expression of the function $\alpha(\omega_i)$ given in \eqref{eq:def_alpha}, we have $\left|\alpha(\omega_i)\right| < 1.03$.
This implies
$$ \left|\alpha(\omega_1)\mu_1 \sin\varphi_1 + \alpha(\omega_2)\mu_2 \sin\varphi_2\right| < 1.03(\left| \mu_1\right| + \left|\mu_2\right|).$$
Therefore, if
\begin{equation*}
\left| \mu_1\right| + \left|\mu_2\right| < 1/1.03,
\end{equation*}
the two surfaces of the crests $\mathcal{C}(I)$ are horizontal for any value of $I_1$ and $I_2$.

%=============================== figure: horizontal crest
%==========================================================

If condition \eqref{eq:horizontal_cond} is not satisfied, $s$ cannot be written as a global function of $\varphi_1$ and $\varphi_2$ in \eqref{eq:crests_cms}, and then two more possibilities arise: b) we can write $\varphi_i$ as a global function of $\varphi_j$ and $s$, or c) we cannot write any variable $\varphi_1$, $\varphi_2$ and $s$ as a global function of the other two and, therefore, the projection of the crest $\mathcal{C}(I)$ on each plane $(\varphi_1 ,\varphi_2)$
, $(\varphi_1 , s)$ and $(\varphi_2 , s)$, has ``holes''.

\begin{figure}[h]
\centering
\subfigure[\emph{Vertical Crest} $\mathcal{C}(I)$ : $\mu_1 =1.7$, $\mu_2 = 0.4$ and $\omega_1 = \omega_2 = 1$.We plotted $\varphi_1,\varphi_2,s\in \left[-\pi/2, 3\pi/2\right)$ for a better illustration. ]
{\includegraphics[scale=0.3]{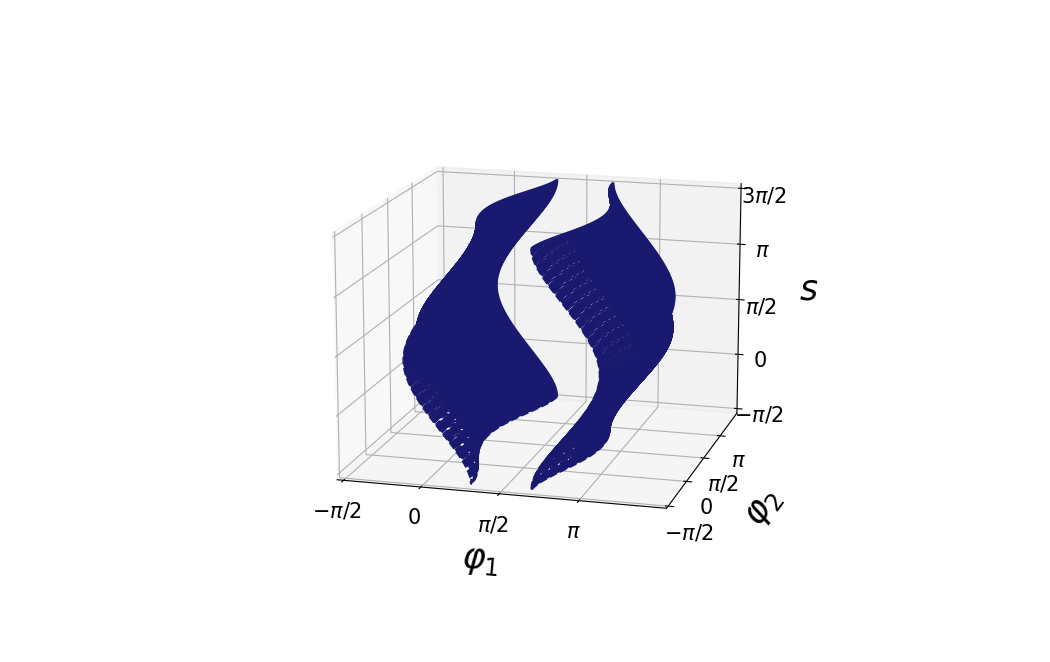}
\label{fig:vertical_crests}
}
\subfigure[\emph{Unseparated Crest} $\mathcal{C}(I)$ : $\mu_1 =0.8$, $\mu_2 = 0.4$ and $\omega_1 = \omega_2 = 1$.We plotted $\varphi_1,\varphi_2,s\in \left[-\pi/2, 3\pi/2\right)$ for a better illustration. ]
{\includegraphics[scale=0.3]{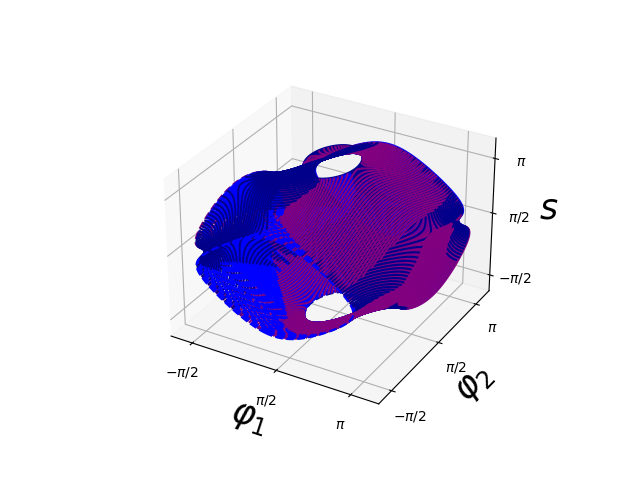}
\label{fig:cristas_nao_separadas}}
\caption{Different kinds of Crest}
\end{figure}

Case b) is only possible if
\begin{equation}
\left|\frac{\alpha(\omega_j)\mu_j}{\alpha(\omega_i)\mu_i}\sin\varphi_j + \frac{\sin s }{\alpha(\omega_i)\mu_i}\right|\leq 1,\label{eq:case_b_vertical}
\end{equation}
for $i, j = 1,2$ and $i\neq j$.
Then, the crest $\mathcal{C}(I)$ is formed by two vertical surfaces $\mathcal{C}_{\text{M}}(I)$ and $\mathcal{C}_{\text{m}}(I)$ that can be parameterized by
\begin{align*}
\varphi_i = \left\{\begin{matrix}
&\eta_{\M}(I, \varphi_j ,s ) :=\arcsin\left( \frac{1}{\alpha(\omega_i)\mu_i}\left( \sin s \right.\right.\\
&-\left.\left. \alpha(\omega_j)\mu_j\sin\varphi_j\right)\right)\mod 2\pi\vspace{0.3cm}\\
&\eta_{\m}(I,\varphi_j , s) :=-\arcsin\left( \frac{1}{\alpha(\omega_i)\mu_i}\left( \sin s\right.\right.\\
&\left.\left. - \alpha(\omega_j)\mu_j\sin\varphi_j\right)\right) + \pi \mod 2\pi.
\end{matrix}\right.
\end{align*}
In this case, $\mathcal{C}(I)$ is called a \emph{vertical} crest; see Fig.~\ref{fig:vertical_crests}.

Case c) only occurs if Eq.~\eqref{eq:horizontal_cond} and \eqref{eq:case_b_vertical} do not hold. Then the crest $\mathcal{C}(I)$ is given by a unique surface; see Fig.~\ref{fig:cristas_nao_separadas}.
Note that the horizontal and vertical crests $\mathcal{C}(I)$ are formed by two disjoint surfaces that can be parameterized separately $\mathcal{C}_{\text{M}}(I)$ and $\mathcal{C}_{\text{m}}(I)$, and $\mathcal{C}(I) = \mathcal{C}_{M}(I)\cup\mathcal{C}_{\text{m}}(I)$.
In case c), $\mathcal{C}(I)$ is called an\emph{unseparated} crest.

%=================== figure: vertical crest
%=======================================================

It is worth noting that to express $\varphi_i$ as a global function of $\varphi_j$ and $s$, $\left|\mu_i\right|$ is needed to be greater than $1/1.03$.
Indeed, assume that there exists an $I$ such that $\varphi_i:=\varphi_i(\varphi_j ,s)$.
From condition \eqref{eq:def_crests}, we have
\begin{equation*}
\sin\varphi_i = -\left(\frac{A_j\omega_j\sin\varphi_j}{A_i\omega_i} + \frac{A_3\sin s}{A_i\omega_i}\right),
\end{equation*}
for any $(\varphi_j , s)\in\mathbb{T}^2$ satisfying \eqref{eq:case_b_vertical}.
In particular, for $\varphi_j = 0 $ and $s = \pm\pi/2$, so
\begin{equation*}
\left| \frac{A_3}{A_i\omega_i}\right| \leq 1, \quad \text{or equivalentely}\quad \left|\frac{1}{\alpha(\omega_i)\mu_i}\right|\leq1.
\end{equation*}
Therefore, we obtain the following.
\begin{equation*}
 \frac{1}{1.03} <\left|\frac{1}{\alpha(\omega_i)}\right| \leq \left|\mu_i\right|.
\end{equation*}
As a consequence, if $\left|\mu_1\right| + \left|\mu_2\right| > 1/1.03$, but $\left|\mu_1\right| , \left|\mu_2\right| < 1/1.03$, there are no vertical crests, that is, there are only horizontal or unseparated crests.

We now summarize all these calculations.
\begin{proposition}
For a fixed value of $I$, if condition \eqref{eq:horizontal_cond} holds, $\mathcal{C}(I)$ is a horizontal crest formed by two disjoint horizontal global surfaces.
Otherwise, we have two cases:
\begin{itemize}
\item[a)] If $|\mu_1| > 1/1.03$ or $|\mu_2|>1/1.03$, $\mathcal{C}(I)$ is a vertical crest formed by two disjoint vertical global surfaces.
\item[b)] If $|\mu_1| + |\mu_2| > 1/1.03$ and $|\mu_1|,\,|\mu_2| < 1/1.03$, $\mathcal{C}(I)$ is an unseparated crest.
\end{itemize}
Moreover, if $|\mu_1| +| \mu_2| < 1/1.03$, then \eqref{eq:horizontal_cond} is satisfied and $\mathcal{C}(I)$ is a horizontal crest formed by two disjoint horizontal global surfaces.
\end{proposition}

\begin{remark}
The values of $\mu_1,\,\mu_2$ providing equalities are bifurcation values for which the two surfaces intersect tangentially.
\end{remark}

\begin{remark}
The notation $\mathcal{C}_{\text{M}}(I)$ and $\mathcal{C}_{\text{m}}(I)$ come from the fact that $(0,0,0)\in \mathcal{C}_{\text{M}}(I)$ and $(\pi,\pi,\pi)\in\mathcal{C}_{\text{m}}(I)$,
and $(0,0,0)$ and $(\pi,\pi,\pi)$ are, respectively, a maximum and a minimum point of the Melnikov potential $\mathcal{L}(I,\varphi,s)$ given by \eqref{eq:mel_pot_gen}, for $a_1,\,a_2,\,a_3>0$.
\end{remark}

\subsubsection{Tangency condition}

The tangency between the NHIM lines $R(I,\varphi,s)$ and the crest $\mathcal{C}(I)$ is an obstacle to the existence of a global scattering map \cite{Delshams2017,Delshams2018}.
Since we only deal with global scattering maps in this paper, we need to avoid such tangencies.
We now make a study about the conditions of their existence.

The crests $\mathcal{C}(I)$ form a family of surfaces, so there exists a tangency between $\mathcal{C}(I)$ and $R(I,\varphi,s)$ if a tangent vector of the straight line $R(I, \varphi,s)$ lies on the tangent bundle of one of these surfaces.

A tangent vector of $R(I,\varphi,s)$ at any point is $-\tilde{\omega}$.
Consider the function $F_I:\mathbb{T}^3\mapsto\mathbb{R}$,
$$F_I(\varphi ,s) = \alpha(\omega_1)\mu_1\sin\varphi_1 + \alpha(\omega_2)\mu_2\sin\varphi_2 + \sin s.$$
We note that the crest $\mathcal{C}(I)$ is defined from \eqref{eq:crests_cms} as the set of $(\varphi,s)\in\mathbb{T}^3$ such that $F_I(\varphi,s) = 0$.
Fixing a point $\mathbf{\Psi} = (\varphi,s)$ in $\mathcal{C}(I)$, the normal vector of $\mathcal{C}(I)$ at the point $\mathbf{\Psi}$ is
\begin{equation*}
\nabla F_I(\mathbf{\Psi}) = (\alpha(\omega_1)\mu_1\cos\varphi_1, \alpha(\omega_2)\mu_2\cos\varphi_2,\cos s).
\end{equation*}
The tangent vector $-\tilde{\omega}$ lies on the tangent space of the crest at the point $\mathbf{\Psi}$ if and only if $\nabla F(\mathbf{\Psi})\cdot \tilde{\omega} = 0$.
This condition is equivalent to
\begin{equation}\label{eq:tang_condition}
\alpha(\omega_1)\omega_1\mu_1\cos\varphi_1 + \alpha(\omega_2)\omega_2\mu_2\cos\varphi_2 + \cos s = 0.
\end{equation}

From \eqref{eq:crests_cms} and \eqref{eq:tang_condition}, there is a tangency between a horizontal crest $\mathcal{C}(I)$ and a \NH line $R(I ,\varphi ,s)$ for $\varphi_1$ and $\varphi_2$ satisfying
\begin{align*}
\left(\sum\limits_{k=1}^2\omega_k\alpha(\omega_k)\mu_k\cos\varphi_k \right)^2
 + \left(\sum\limits_{k=1}^2\alpha(\omega_k)\mu_k\sin\varphi_k \right)^2 = 1.
\end{align*}

Denote by
\begin{align*}
f_{I}(\varphi) = \left(\omega_1\alpha(\omega_1)\mu_1\cos\varphi_1 + \omega_2\alpha(\omega_2)\mu_2\cos\varphi_2\right)^2\\ + \left(\alpha(\omega_1)\mu_1\sin\varphi_1 + \alpha(\omega_2)\mu_2\sin\varphi_2\right)^2.
\end{align*}

\begin{figure}[h]
\begin{center}
\includegraphics[scale=0.28]{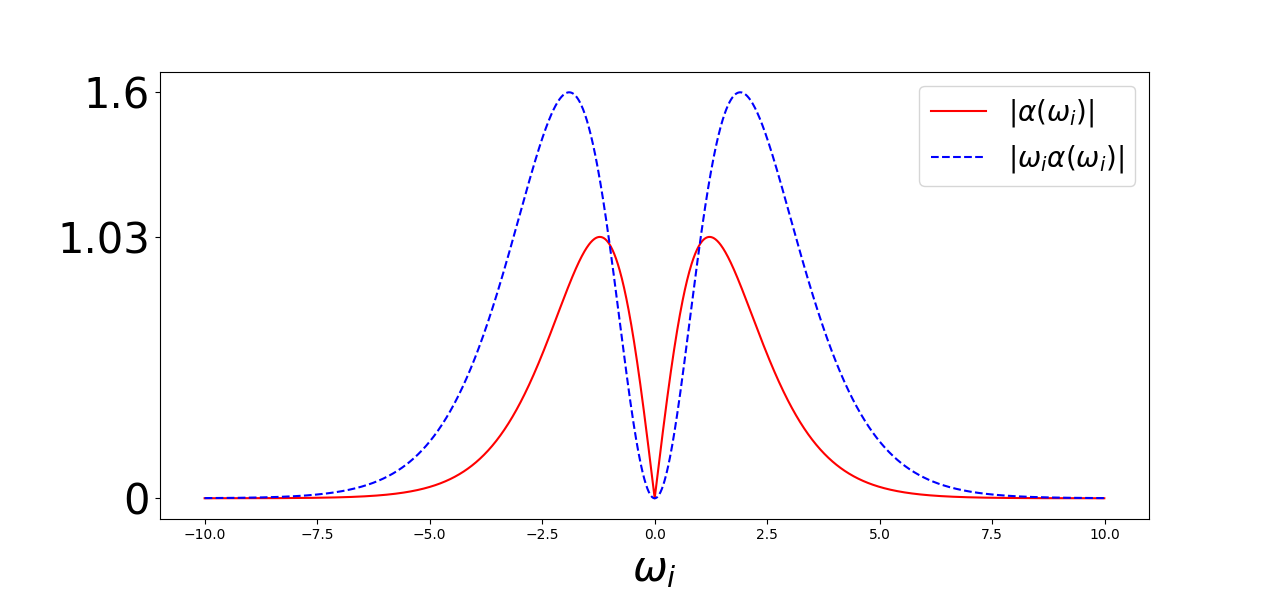}
\end{center}
\caption{Graph of $\left|\alpha\right|$ and $\left|\omega_i\alpha\right|$.}
\label{fig:alpha_and_beta}
\end{figure}

Note that for values of $\mu_1$ and $\mu_2$ such that $f_I(\varphi)< 1$ there is no tangency.
As commented above, $\left|\alpha(\omega_i)\right| < 1. 03$.
From \eqref{eq:def_alpha}, we obtain $\left|\omega_i\alpha(\omega_i)\right| < 1.6$, see Fig.~\ref{fig:alpha_and_beta}.
This implies
\begin{align*}
f_I(\varphi)&<(1.6)^2\left(\left|\mu_1\cos\varphi_1\right| + \left|\mu_2\cos\varphi_2\right|\right)^2 \\&+ (1.03)^2\left(\left|\mu_1\sin\varphi_1 \right|+ \left|\mu_2\sin\varphi_2\right|\right)^2 \\
 &<1.6^2 \left|\mu_1\right|^2 + 1.6^2\left| \mu_2\right|^2\\& + 2\left|\mu_1\right|\left|\mu_2\right|\left(1.6^2\left|\cos\varphi_1\right|\left|\cos\varphi_2\right|\right.\\&\left.  + 1.03^2\left|\sin\varphi_1\right|\left|\sin\varphi_2\right|\right)\\
 &< 1.6^2(\left|\mu_1\right| + \left|\mu_2\right|)^2.
\end{align*}
Therefore it is enough to require $\left|\mu_1\right| + \left|\mu_2\right| < 1/1.6 = 0.625$ to ensure $f_I(\varphi)<1$ for any value of $I$.
It is easy to verify that if $\left|\mu_1\right| + \left|\mu_2\right| > 1/1.6 $ there exist $I$ and $\varphi$ such that $f_I(\varphi)>1$.

We now collect all these properties in the following proposition.

\begin{proposition}\label{prop:bifur_crests}
For any I, consider the crest $\mathcal{C}(I)$ defined in \eqref{def:crest} and the \NH lines $R(I,\varphi ,s)$ defined in \eqref{eq:nhim_line}.
\begin{itemize}
\item[a)] For $\left|\mu_1\right| + \left|\mu_2\right| <0.625$ the crest $\mathcal{C}(I)$ is formed by two horizontal global surfaces and the intersections between any surface and any \NH line are transversal.
\item[b)] For $0.625 \leq \left|\mu_1\right| + \left|\mu_2\right|\leq 1/1.03$ the crest $\mathcal{C}(I)$ is still formed by two horizontal global surfaces, but for some values of $I$ there are \NH lines that are tangent to the surfaces.
\item[c)] For $ 1/1.03 < \left|\mu_1\right| + \left|\mu_2\right|$ and $\left|\mu_1\right|,\left|\mu_2\right| \leq 1/1.03$, the crest $\mathcal{C}(I)$ is either formed by two global horizontal surfaces or is unseparated, and for some values of $I$ there are \NH lines that are tangent to the crest.
\item[d)]For $1/1.03 < \left|\mu_1\right|$ or $1/1.03 < \left|\mu_2\right|$ , the crest $\mathcal{C}(I)$ is either formed by two horizontal global surfaces, either formed by two vertical global ones, or it is unseparated, and for some values of $I$ there are \NH lines that are tangent to the crest.
\end{itemize}
\end{proposition}

Throughout this paper, we will consider only case a), that is, $\left|\mu_1\right| + \left|\mu_2\right| <0.625$ which is contained in hypothesis \eqref{eq:parametros} of Theorem 1, because in this way the crest is formed by two surfaces and the \NH lines are transverse to them.

We have assumed until now that $s\in \mathbb{T}$. This implies that $\mathcal{C}(I) = \mathcal{C}_{\text{M}}(I)\cup\mathcal{C}_{\text{m}}(I)$ and $R(I,\varphi,s)$ intersects infinitely many times $\mathcal{C}_{\text{M}}(I)$ and $\mathcal{C}_{\text{m}}(I)$.
We will restrict our analysis to a finite number of intersections to obtain our results.
Taking $s\in\mathbb{R}$, the crest $\mathcal{C}(I)$ is no longer formed by two horizontal surfaces, but it is an infinitely countable family of horizontal surfaces, that is, $\mathcal{C}(I) = \bigcup_{j\in\mathbb{Z}}\mathcal{C}_j(I)$, where $\mathcal{C}_j(I)$ can be parametrized by
\begin{align}\label{eq:par_j_even}
\xi_j(I,\varphi) =-\arcsin\left(\sum\limits_{k=1}^2\mu_k\alpha(\omega_k)\sin\varphi_k\right) + 2\pi j,
\end{align}
if $j$ is even, otherwise
\begin{align}\label{eq:par_j_odd}
\xi_j(I,\varphi) = \arcsin\left(\sum\limits_{k=1}^2\mu_k\alpha(\omega_k)\sin\varphi_k\right)+\pi + 2\pi j.
\end{align}

\begin{remark}
Note that if $s\in\mathbb{T}$, $\mathcal{C}_j\equiv \mathcal{C}_{\text{M}}$ for all even $j$.
Analogously, $\mathcal{C}_j\equiv \mathcal{C}_{\text{m}}$ for all odd $j$.
\end{remark}

Therefore, we immediately have the following corollary from Prop.~\ref{prop:bifur_crests}:
\begin{corollary}
\label{cor:horizontal_crests}
Consider the crest $\mathcal{C}(I)$ defined in \eqref{def:crest} and the \NH lines $R(I,\varphi ,s)$ defined in \eqref{eq:nhim_line} for $s\in\mathbb{R}$.
For $\left|\mu_1\right| + \left|\mu_2\right| <0.625$ the crest $\mathcal{C}(I)$ is formed by infinitely countable horizontal surfaces $\mathcal{C}_j(I)$, $j\in\mathbb{Z}$, that is, $\mathcal{C}(I)= \bigcup_{j\in\mathbb{Z}}\mathcal{C}_j(I)$, and the intersections between any surface $\mathcal{C}_j(I)$ and any \NH line $R(I,\varphi,s)$ are transversal.
\end{corollary}

After this construction, each line $R(I,\varphi, s)$ intersects each $\mathcal{C}_j(I)$ at a unique point.
We denote by $\tau^*_j$ the value of $\tau$ such that $R(I,\varphi,s)$ intersects $\mathcal{C}_j(I)$.
Following the same philosophy, we denote $\mathcal{L}_j^*(I,\theta) = \mathcal{L}(I,\varphi-\tau^*_j\omega, s-\tau^*_j)$ and $\mathcal{S}_j$ the scattering map associated to $\mathcal{C}_j$, that is,
\begin{equation*}
\mathcal{S}_j(I,\theta)^{\text{T}} = (I,\theta)^{\text{T}} + \varepsilon J\nabla \mathcal{L}_j^*(I,\theta)^{\text{T}} + \mathcal{O}(\varepsilon^2),
\end{equation*}
where $J=\left(\begin{smallmatrix}0&\text{Id}_2\\-\text{Id}_2&0\end{smallmatrix}\right)$ is the $4\times4$ symplectic matrix.

\begin{corollary}\label{cor:global_scattering}
Under the same conditions of Corollary~\ref{cor:horizontal_crests}, the scattering map $\mathcal{S}_j(I,\theta)$ associated to $\mathcal{C}_j$ is globally defined, i.e., it is well defined for any $(I,\theta)\in \mathbb{R}^2\times \mathbb{T}^2$.
\end{corollary}

\begin{proof}
The existence of $\tau_j^*(I,\theta)$ for all $(I,\theta)\in\mathbb{R}^2\times \mathbb{T}^2$ holds by the fact that the slope of $R(I,\varphi,s)$ is $\tilde{\omega}~=~(\omega_1,\omega_2,1)$ and the surfaces $\mathcal{C}_j(I)$ are horizontal and there is no tangency between $R(I,\varphi,s)$ and $\mathcal{C}_j(I)$.
\end{proof}

\begin{remark}
In practice, it is enough to deal with $\mathcal{S}_0$ and $\mathcal{S}_1$, which are also called primary scattering maps\cite{Delshams2017,Delshams2018}.
Indeed, in this paper we will only explicitly use~$\mathcal{S}_0$.
\end{remark}

\subsection{ Transversality between the Inner and scattering flows}

To perform diffusion through a scattering map, that is, to find trajectories of Hamiltonian system \eqref{eq:ham}$+$\eqref{eq:pert} where the action $I$ changes significantly, we would like to prove that an orbit of the scattering map is not confined in an orbit of the inner dynamics, since the value of $I$ is almost constant under the action of the inner map.

Recall that the scattering map is the $\varepsilon$-time flow of the Hamiltonian $-\mathcal{L}_j^*(I, \theta)$ up to order $\mathcal{O}(\varepsilon^2)$.
Then, from a geometrical point of view, looking at the flow associated to the reduced Poincar\'e function $\mathcal{L}_j^*$ will be very useful.

\begin{figure}
\begin{center}
\includegraphics[scale = 0.45]{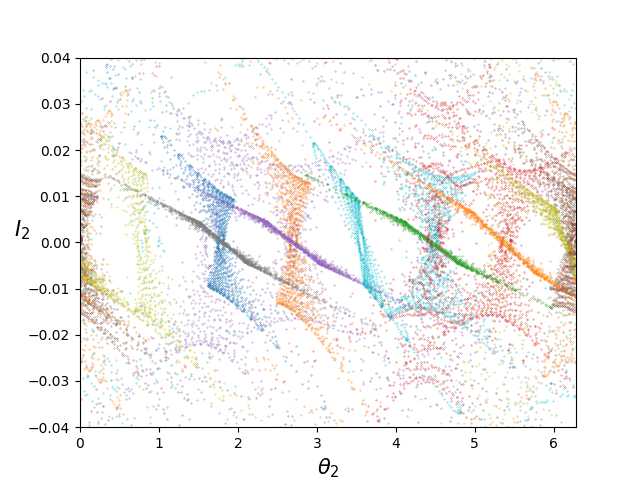}
\end{center}
\caption{Dynamics of the Poincar\'e map the section
$I_1=0$ at the energy level $\mathcal{L}^*(0 , 0 , 5 \pi/4 , 5\pi / 4)$ . Initial conditions are taken at the section
$\left\{ \left( I_1,I_2,\theta_2
\right)=(0, 0,[5\pi / 4 - 0.2 , 5\pi / 4 +0.2  ]) \right\}$ for $\mu_1 = 0.2$ and $\mu_2 = 0.3$. Each color represents a different orbit of the Poincar\'e map.}
\label{fig:non_int}
\end{figure}

\begin{definition}
The flow $\phi_t^{-\mathcal{L}_j^*}$ of the Hamiltonian $-\mathcal{L}_j^*$ given by the equations
\begin{equation}\label{eq:l^*_ham_equations}
\dot{I} = \frac{\partial \mathcal{L}_j^*}{\partial \theta}(I,\theta)\quad\quad\text{and}\quad\quad \dot{\theta} = -\frac{\partial \mathcal{L}_j^*}{\partial I}(I, \theta)
\end{equation}
is called the \emph{scattering} flow associated to $\mathcal{L}_j^*$.
\end{definition}

The scattering flow is 4-dimensional and autonomous
and exhibits rich dynamics. In Figure~\ref{fig:non_int} we show a
simulation of the Poincar\'e map onto a hyperplane given by
$I_1=\text{constant}$, where non-integrability is clearly manifested.

From \eqref{eq:reduced_poincare_explicit}, the above Hamiltonian equations take the form, for $i = 1,2$,
\begin{align}\label{eq:edo_SM_k=0}
\dot{I}_i &= -A_i \sin(\theta_i - \omega_i\tau_j^*)
\\
\dot{\theta}_i &= -\Omega_i\left(\frac{dA_i}{d\omega_i}\cos(\theta_i-\omega_i\tau_j^*) + \tau^*A_i\sin(\theta_i-\omega_i
 \tau_j^*) \right).\nonumber
\end{align}

\begin{remark}\label{rem:error}
Take an initial point $(I_0,\theta_0)$. After a large number $k$ of iterates of scattering map $\mathcal{S}_j$, $(I_k,\theta_k)=\mathcal{S}^k_j(I_0,\theta_0)$ is expected to accumulate an error with respect to the scattering flow $\phi_{t}(I_0,\theta)$ at time $t = k\varepsilon$, as the flow is just an approximation of the scattering orbit. However, this error can be controlled by also using the inner dynamics in order to keep the iterations of $\mathcal{S}_j$ close to the scattering flow.
A complete description of this mechanism and an analysis of this accumulated was done for 2$+$1/2 d.o.f.\citep{Delshams2017}
\end{remark}

\begin{lemma}\label{lem:fixed_points}
The scattering flow with equations \eqref{eq:edo_SM_k=0} has only four equilibrium points: $(0,0,0,0)$, $(0,0,0,\pi)$, $(0,0,\pi,0)$ and $(0,0,\pi,\pi)$.
\end{lemma}
\begin{proof}
It is an immediate consequence of \eqref{eq:edo_SM_k=0} and the form of $A_i$, given in \eqref{eq:coef_A}, and their derivatives.
\end{proof}

We can now state the following result.

\begin{proposition} \label{prop:scattering_inner_transversal}
Consider Hamiltonian~\eqref{eq:ham}$+$\eqref{eq:pert}.
For $\left|\mu_1\right| + \left|\mu_2\right| < 0.625$ and $I$ not $\varepsilon$-close to $(0,0)$, there is no common orbit of the inner flow and the scattering flow.
%a scattering flow associated to $-\mathcal{L}_j^*$, with $|\tau^*_j| \ll 1/\varepsilon $, does not contain any orbit of the inner map.
\end{proposition}

\begin{proof}
From \eqref{eq:inner_eq}, the first integrals on variables $(I_i,\theta_i)$ are given by
\begin{equation}\label{eq:first_integral}
F_i = \frac{\Omega_i I_i^2}{2} + \varepsilon a_i\cos\theta_i,\quad i = \left\{1,2\right\}.
\end{equation}
The transversality of any invariant set of the inner flow and the scattering flow holds if the gradient of the level surfaces of $F_1$ and $F_2$ are not parallel to the gradient of the level surfaces of $\mathcal{L}^*$, or equivalently,
\begin{align}\label{eq:poisson_bracket}
\left\{F_i, -\mathcal{L}_j^*\right\}(\phi^{-\mathcal{L}_j^*}_t(I^0,\theta^0)) \neq 0,\quad i = \left\{1,2\right\},
\end{align}
where $\left\{,\right\}$ is the Poisson bracket and $\phi^{-\mathcal{L}_j^*}_t(I^0,\theta^0)$ is the scattering flow associated to $-\mathcal{L}_j^*$.
Note that
\begin{align*}
&\left\{F_i , -\mathcal{L}_j^*\right\} =- \frac{\partial F_i}{\partial \theta_i}\frac{\partial \mathcal{L}_j^*}{\partial I_i} + \frac{\partial F_i }{\partial I_i}\frac{\partial \mathcal{L}_j^*}{\partial \theta_i}\\
&= \varepsilon a_i \Omega_i \sin\theta_i \left(\frac{dA_i}{d\omega_i}\cos(\theta_i - \omega_i\tau_j^* ) + \tau_j^*A_i\sin(\theta_i - \omega_i\tau_j^*)\right)\\
& - \omega_i A_i \sin(\theta_i- \omega_i\tau_j^*)
\end{align*}
Suppose that \eqref{eq:poisson_bracket} does not hold.

In the case that $\left|I_i^0\right|\gg \varepsilon$, for $i=1,2$, the dominant part is $-\omega_iA_i\sin(\theta_i-\omega_i \tau^*)$.
So, $\left\{F_i , -\mathcal{L}^*\right\} = 0$ only if $\sin(\theta_i -\omega_i\tau^*) = 0$, for $i= 1,2$.
This implies that, from \eqref{eq:edo_SM_k=0}, $I_i$ is constant, for $i= 1,2$.
As $F_i$ and $I_i$ are constant, $\theta_i$ is also constant, $i = 1,2$.
Then, we can conclude that $(I,\theta)$ is an equilibrium point of $\mathcal{L}_j^*$.
From Lemma~\ref{lem:fixed_points}, $I = (0,0)$, a contradiction.

Now, consider the case that $\left|I_i\right|\gg \varepsilon$ and $I_l \approx \varepsilon$ , for $i\neq l$.
As we have seen before, $I_i$ and $\theta_i$ are constant.
This implies that we can reduce the scattering flow two to a 2D flow.
The transversality relies on the same argument for the 2$+$1/2 d.o.f. case.\cite{Delshams2011,Delshams2017}

\end{proof}
\begin{corollary}
The same proposition is valid for the scattering map, since the scattering map is the scattering flow plus a term of order $\mathcal{O}(\varepsilon^2)$ for $\varepsilon$ small enough.
\end{corollary}
%------------------------------------------------------
% DIFFUSION
%---------------------------------------------------
\subsection{Construction of diffusing paths}

We are interested in finding a finite drift in the 2-dimensional variable $I$.
In general, we expect to obtain an increment in the value of $I$ by iterating the scattering map.
Therefore, the points where the action variable $I$ does not change under the action of the scattering map are useless for our mechanism.
From Lemma~\ref{lem:fixed_points} and Prop.~\ref{prop:scattering_inner_transversal}, there are four points to avoid: $(0,0,0,0),\,(0,0,\pi,0),\,(0,0,0,\pi)$ and $(0,0,\pi,\pi)$ since they are equilibrium points of all scattering flows and inner dynamics.

As we are assuming $\mu_1$, $\mu_2$ such that $\left|\mu_1\right| + \left|\mu_2\right|< 0.625$, we find that all surfaces $C_j(I)$ are horizontal for any value of $I$ and are transversally intersected by $R(I,\varphi,s)$, see Corollary~\ref{cor:horizontal_crests}.
Therefore, $\tau_j^*(I,\theta)$ is always defined, see the proof of Corollary~\ref{cor:global_scattering}, and we can introduce a new variable:
\begin{equation}\label{def:psi}
\psi_j:= \theta - \tau_j^*(I,\theta)\omega.
\end{equation}
By \eqref{eq:crests_cms},\eqref{eq:par_j_even} and \eqref{eq:par_j_odd}
\begin{align*}
s-\tau_j^*(I\varphi,s) = \xi_j(I,\varphi-\tau_j^*(I,\varphi,s)\omega)
\end{align*}
By \eqref{eq:def_theta} and \eqref{eq:tau_theta}
\begin{align*}
-\tau_j^*(I,\theta) = \xi_j(I,\theta-\tau_j^*(I,\theta)\omega).
\end{align*}
This immeadiately implies that $\psi$ has a well-defined inverse function
\begin{equation*}
\theta = \psi -\xi_j(I,\psi_j)\omega.
\end{equation*}
This variable will be helpful in the proof of the following theorem.

%
%-----------------------
\section{Global Instability}\label{sec:global_instability}

\begin{reptheorem}{thm:main_theorem}
Consider the Hamiltonian \eqref{eq:ham}$+$\eqref{eq:pert}.
Assume $a_1a_2a_3 \neq 0$ and $\left|a_1/a_3\right| + \left|a_2/a_3\right| < 0.625$.
Then, for every $0<\delta < 1$ and $R >0$,  there exists $\varepsilon_0=\varepsilon_0(\delta,R)>0$ such that for every $I_+, I_-$ satisfying $\left|I_{\pm}\right|< R$,
there exists an orbit $\tilde{x}(t)$ and $T >0$, such that
\begin{align*}
\left|I(\tilde{x}(0)) - I_- \right| &\leq C\delta\\
\left|I(\tilde{x}(T)) - I_+ \right| &\leq C\delta \nonumber
\end{align*}
\end{reptheorem}

\begin{proof}

Actually, we are going to prove that given $I_-$, $I_+$ and a curve $\gamma:\left[0,t_\text{f}\right]\rightarrow \mathbb{R}^2$, such that $\gamma(0) = I_-$ and $\gamma(t_\text{f}) = I_+$, there exists an orbit $\tilde{x}(t)$ of the Hamiltonian system given by \eqref{eq:ham}$+$\eqref{eq:pert} satisfying
\begin{equation}\label{eq:theo_1_proof}
\left|I\left(\tilde{x}(t)\right) - \gamma(\alpha(t))\right|\leq C\delta,\quad t\in\left[0,t_\text{f}\right].
\end{equation}

Our strategy will be to use iterations of the inner and scattering maps to construct a pseudo-orbit $\delta$-close to $\gamma$.
Later, we will apply shadowing lemmas to ensure the existence of a real diffusing orbit in which condition \eqref{eq:theo_1_proof} is satisfied.

Consider first the case that $\gamma$ is a horizontal segment in the plane $(I_1,I_2)$, so that $I_{-2} = I_{+2}$,
where $v_2=\partial_{\theta_2}\mathcal{L}^*(I, \theta)$ does not vanish.
Given $\delta <l(\gamma)$, where $l(\gamma)$ is the length of $\gamma$, we define $N:= \left \lfloor \frac{l(\gamma)}{\delta} \right \rfloor$ and take the finite open cover of the image of the curve $\gamma$ given by
\begin{equation*}
\bigcup_{i = 0}^{N}B_{\delta}(\gamma(t_i)),
\end{equation*}
where $B_{\delta}(\gamma(t_i)) = \left\{p \in \mathbb{R}^2 : \left\|\gamma(t_i)- p\right\|_{\infty} < \delta\right\}$ and $t_0 = 0$.

Let $I^i\in B_{\delta}(\gamma(t_i))\setminus B_{\delta}(\gamma(t_{i+1}))$.
Applying the scattering map, we want to obtain a $I^{j}\in B_{\delta}(\gamma(t_{i+1}))$.
From now on, we denote by $I^k$ the values of $I$ obtained by the $k$th iteration under the scattering
map $\mathcal{S}$. We should ensure $I^k \in B_{\delta}(\gamma(t_i)) \cup B_{\delta}(\gamma(t_{i+1}))$
for $k \in \left\{1, \dots , j\right\}.$
Recall that the scattering map can be written as
\begin{equation*}
\mathcal{S}(I,\theta)^{\text{T}} = (I,\theta)^{\text{T}} + \varepsilon J\nabla \mathcal{L}^*(I,\theta)^{\text{T}} + \mathcal{O}(\varepsilon^2),
\end{equation*}
where $J$ is the symplectic matrix, and the term $\mathcal{O}(\varepsilon^2)$ can be bounded by $\varepsilon^2 M$, for some constant $M$ depending just on $R$.
Therefore, defining the vector $u^i = \gamma(t_{i+1})- I^i$, we have to obtain a $\theta^*$ such that
\begin{equation}\label{eq:condition_direction}
u^i_j \partial_{\theta_j}\mathcal{L}^*(I^i, \theta^*) >0, \text{ for } j = \left\{ 1,\,2\right\},
\end{equation}
where $\partial_{\theta_j}\mathcal{L}^*(I, \theta) = -A_j\sin\left(\theta_j - \omega_j\tau^*(I,\theta)\right)$ with $A_j$ defined in \eqref{eq:coef_A}, see Fig.~\ref{fig:balls}.
Assuming $a_1,\,a_2>0$, \eqref{eq:condition_direction} is satisfied if and only if $\theta^* -\tau^*(I^i,\theta^*)\omega^i \in \left(\pi, 2\pi\right)^2$.

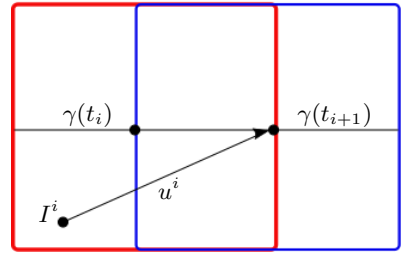
\begin{figure}[h]
\centering
\begin{overpic}[scale=0.45,unit=1mm]{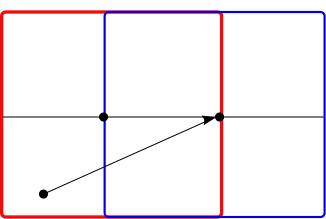}
\put(-5,18){{\parbox{0.35\linewidth}{%
     \begin{equation*}
     \gamma(t_i)
     \end{equation*}}}}
\put(30,18){{\parbox{0.3\linewidth}{%
     \begin{equation*}
     \gamma(t_{i+1})
     \end{equation*}}}}
 \put(-8,5){{\parbox{0.3\linewidth}{%
     \begin{equation*}
     I^i
     \end{equation*}}}}
 \put(8,8){{\parbox{0.3\linewidth}{%
     \begin{equation*}
     u^i
     \end{equation*}}}}
\end{overpic}
\caption{$\partial B_{\delta}(\gamma(t_i))$ in red and $\partial B_{\delta}(\gamma(t_{i+1}))$ in blue.}
\label{fig:balls}
\end{figure}

If for initial values $(I^i, \theta^i)$, $\theta^i - \tau^*(I^i,\theta^i)\omega_i \notin (\pi,2\pi)^2$, we need to use inner dynamics to displace $\theta^i$ to a point $\theta^*$ such that
\begin{equation*}
\psi^*: = \theta^* - \tau^*(I^i,\tau^*)\omega^i \in (\pi, 2\pi)^2.
\end{equation*}
The inner dynamics is very simple, and we can assume that it is approximately horizontal for finite time\cite{Delshams2017}, i.e., it is described by the equations
\begin{equation*}
\dot{I}_j = 0 \quad\text{and}\quad\dot{\varphi}_j = \omega_j,\quad j = 1,2.
\end{equation*}
Therefore, $\varphi(t)  = \omega t + \varphi(0)$.
Then, the existence of $\theta^*$ is equivalent to the existence of $t^*$ satisfying
\begin{equation*}
\theta(t^*) - \tau^*(I^i ,\theta(t^*))\omega^i\in (\pi , 2\pi)^2,
\end{equation*}
where $\theta(t) = \theta^i + t\omega^i $.
Define
$$\psi(t) = \theta^i(t) - \tau^*(I^i , \theta^i(t))\omega^i.$$

Assume $\omega^i_1 \geq \omega^i_2$ without loss of generality, and as we can choose the initial point with a $\delta$-error, we can take $I^i_1,\,I^i_2\neq 0$.
We can write
\begin{equation}\label{eq:psi_2 as_fctin_psi_1}
\psi_2 = \frac{\omega^i_2}{\omega^i_1}\psi_1 + \psi_0,\quad\quad \psi_0 := \theta^i_2 - \omega^i_2\theta^i_1/\omega^i_1.
\end{equation}

For $\omega^i_2/\omega^i_1 \in \mathbb{R}\setminus \mathbb{Q}$, $\left(\psi_1 , \psi_2( \psi_1 )\right)$ is dense in $\mathbb{T}^2$, then there exists a $t^*$ such that $\left(\psi_1(t^*) , \psi_2(t^*)\right)\in \left(\pi , 2\pi\right)^2$.

For $\omega^i_2/\omega^i_1  = p/q \in \mathbb{Q}$, consider $q\,p>0$, the other case is analogous.
This implies $0< p/q \leq 1$.

We fix $\psi_1 = 2\pi l$, $l\in \mathbb{Z}$.
For these values of $\psi_1$, $\psi_2$ in \eqref{eq:psi_2 as_fctin_psi_1} can be rewritten as
\begin{align*}
\psi_2(l) = r_l(\psi_0) = 2\pi\left(\frac{p}{q}\right) l + \psi_0.
\end{align*}
Note that $\psi_2$ is equivalent to a rotation of $\psi_0$ by the angle $2\pi p/q$ on $S^1$.
As $r_{l}(\psi_0)$ is a $q$-periodic function, an orbit $\mathcal{O} = \left\{ \psi_0, \dots , r_{q-1}(\psi_0)\right\}$ is an equidistant cover of $S^1$. Therefore, if $q\neq 1$ or $\psi_0\neq \pi$, there exists $l'\in \left\{0, \dots, q-1\right\}$ such that $r_{l'}(\psi_0)\in (\pi,2\pi)$.
This immediately implies that $\left(\psi_1 , \psi_2(\psi_1)\right)$ intersects $\left(\pi , 2\pi\right)^2$.

For $q = 1$ and $\psi_0 = \pi$, then $\omega^i_1 = \omega^i_2$ and it is easy to verify that $(\psi_1 , \psi_2(\psi_1))$ does not intersect $(\pi , 2\pi]^2$.
Therefore, in this specific case, the inner dynamics is not enough for our purpose.
So, here our strategy is slightly different.
We will use the scattering map in a similar way to the inner map: to change the values of $\theta$ while the value of $I$ remains almost fixed.

First, we apply the inner dynamics up to $\psi_1 = 0$.
From \eqref{eq:psi_2 as_fctin_psi_1}, this implies $\psi_2 = \pi$.
For these values, we have $\dot{I}= 0$.
From Prop.~\ref{prop:scattering_inner_transversal}, when we apply the scattering map on this point, its image is on another torus of the inner dynamics.
As the values of the variable $I$ remain fixed up to $\mathcal{O}(\varepsilon^2)$, we can consider that these values are within $B_{\delta}(\gamma(t_i))$.
Out of this problematic torus, we can apply the algorithm developed previously.

Now we want to prove that $I^k$ is $\delta$-close to the curve $\gamma$ for any $k\in\left\{1,\dots,j\right\}$.
We have
\begin{equation}
I^{k} = I^{k-1} + \varepsilon \partial_{\theta}\mathcal{L}^*(I^{k-1}, \theta^{k-1}) + \mathcal{O}(\varepsilon^2).\label{eq:I_i_+1}
\end{equation}
So, $I^{k}$ is $\delta$-close to $\gamma$ if the following condition is satisfied.
$$
\left|I_2^{k}-\gamma_2(t_{i + 1 })\right| \leq \delta.
$$
If we assume $I^{k-1}_2\leq \gamma_2(t_{i+1})$, we only need to verify $I^{k}_2\leq \gamma_2(t_{i+1}) + \delta$.

From \eqref{eq:I_i_+1} and if we consider only the terms of the first order, the condition is equivalent to
$$ \varepsilon  < \frac{\gamma_2(t_i) - I_2^k +\delta}{v^k_2-\varepsilon M}.$$

Note that $\gamma_2(t_i)-I_2^k < \delta$, then the numerator is positive.
In addition, $v_2^k \geq m_2(\gamma)$, the minimum value of $|v_2|$ along the curve $\gamma$, so it is enough to require
$\varepsilon<m_2(\gamma)/(2M)$, as well as
\begin{equation}
\varepsilon < 2\frac{\gamma_2(t_i)-I_2^k + \delta}{m_2(\gamma)}. \label{eq:epsilon_1}
%\varepsilon < \frac{\gamma_2(t_i)-I_2^k + \delta}{\left\|v^k\right\|_{\infty}}. \label{eq:epsilon_1}
\end{equation}
Define $\varepsilon_i = \sup\left\{\varepsilon : \varepsilon < \frac{\gamma_2(t_i)-I_2^i+ \delta}{m_2(\gamma)} \right\}$, we obtain that for any $0<\varepsilon \leq \varepsilon_i$, $I^{k}$ is $\delta$-close to $\gamma$.

For $u_2 \leq 0$, \eqref{eq:epsilon_1} takes the form
\begin{equation*}
\varepsilon < 2\frac{I_2^k-\gamma_2(t_i) + \delta}{m_2(\gamma)}.
%\varepsilon < \frac{I_2^i -\gamma_2(t_i) + \delta}{\left\|v^k\right\|_{\infty}}.
\end{equation*}

Now we wish to obtain a similar result for any iterate of a scattering map.
Introducing $m=\min_{|I|\leq R}v$, the result is true for $\varepsilon<1/(2M)$ and
\begin{equation*}
\varepsilon < \frac{2\delta}{m}.
\end{equation*}
That is, if we take $\varepsilon_0 = \min\{\frac{1}{2m},\frac{2\delta}{m}\}$, for any $\varepsilon < \varepsilon_0$ we obtain a pseudo-orbit $\delta$-close to $\gamma$.

For vertical lines, the same result can be stated \emph{mutatis mutandis}.

The above argument holds for any $\gamma$ that is $\delta$-distant from the equilibrium proint of scatterig flow, see Lemma \ref{lem:fixed_points}. If such condition it is not satified, it can be approximated by $\gamma'$ that has a distance $\mathcal{O}(\delta/2)$ from $\gamma$ close to those points.

For a more general case, that is, $C^1$-curve $\gamma: \left[0 , t^*\right] \rightarrow \mathbb{R}^2$ such that $\gamma(0)~= I_-$, $\gamma(t^*) = I_+$, we take a stairstep curve $\gamma_{\text{step}}$, a combination of horizontal and vertical lines, in such a way that $\gamma_{\text{step}}$ is a good enough approximation of $\gamma$, where `good enough' we mean, the result holds for $\gamma$ applying the above results (for horizontal and vertical lines) for $\gamma_{\text{step}}$.

We can now apply shadowing techniques well adapted to NHIM~\cite{fontich2000,fontich2003,GideaLS20}, due to the fact that the inner dynamics is simple enough to satisfy the required hypothesis of these references, to ensure the existence of a diffusion trajectory.
\end{proof}

%============================SECTION:HIGHWAYS==========================
\subsection{Highways}\label{sec:Highways_3}

The existence of a special invariant set of $-\mathcal{L}_0^*$, called \emph{Highway}, was very useful for the study of the case with 2$+$1/2 d.o.f.\cite{Delshams2017}.
Iterations of a scattering map along a Highway were enough to obtain a large drift on the action variable $I$.
Moreover, we estimated the time of diffusion of these orbits.
We want to detect a similar invariant set for our 3$+$1/2 model with the same goal in mind.

We define a Highway as an invariant set $\mathcal{H} = \left\{ (I,\Theta(I) ) \right\}$ of the Hamiltonian given by the reduced Poincar\'{e} function $\mathcal{L}_0^*(I,\theta)$ which is contained in the level energy $\mathcal{L}_0^*(I,\theta) = A_3$.
Therefore, it is a Lagrangian manifold, that is, $\Theta(I)$ is a gradient function, so there exists a function $F(I)$ such that $\Theta(I) = \nabla F(I)$.
As $\Theta$ is a gradient function, it must satisfy the following condition.
\begin{equation*}
\frac{\partial \Theta_1}{\partial I_2} = \frac{\partial \Theta_2}{\partial I_1}.
\end{equation*}
This condition is equivalent to
\begin{equation*}
\frac{\partial^2 F}{\partial I_2\partial I_1 } = \frac{\partial^2 F}{\partial I_1 \partial I_2}.
\end{equation*}

\begin{remark}
We take the energy level $\mathcal{L}_0^*(I,\theta) = A_3$ because we have from \eqref{eq:reduced_poincare_explicit} that
\begin{equation*}
\lim_{(I_1,I_2)\rightarrow (\pm \infty , \pm \infty)}\mathcal{L}^*_0(I,\theta) = A_3.
\end{equation*}
In addition, $\mathcal{L}^*_0(0,0,\pm\pi/2, \pm\pi/2) = A_3$.
Therefore, the level surface $\mathcal{L}^*_0(I,\theta) = A_3$ is the unique connection possible of the origin of the plane $(I_1,I_2)$ to arbitrarily large values of $\left|I_1\right|$ and $\left|I_2\right|$.
\end{remark}

We start by proving the local existence of Highways for $I_1,I_2\gg 0$.

\begin{proposition}[Local existence of Highways]
\label{prop:high_asymp}
Consider the Hamiltonian \eqref{eq:ham}$+$\eqref{eq:pert}.
Assume $a_1a_2a_3 \neq 0$ and $\left|a_1/a_3\right| + \left|a_2/a_3\right| < 0.625$.
For $\left|I_1\right|$ and $\left|I_2\right|$ close to infinity, the function $F$ takes the asymptotic form
\begin{align}\label{eq:lagrange_candidate}
F(I) &= \frac{\pi}{2}\left(\pm I_1 \pm I_2\right) -\sum_{i = 1}^2\frac{2a_i \sinh(\pi/2)}{\pi^4\Omega_i}\left( \pi^3\left|\omega_i\right|^3\right.\nonumber\\
&\left. + 6\pi^2\omega_i^2 + 24\pi\left|\omega_i\right| + 48  \right)e^{-\pi\left|\omega_i\right|/2} \nonumber\\
& + \mathcal{O}(\omega_1^2\omega_2^2e^{-\pi(\left|\omega_1\right| + \left|\omega_2\right|)/2}),
\end{align}
\end{proposition}

\begin{proof}
Assume a candidate for a function $F(I)$ given by \eqref{eq:lagrange_candidate},
such that $\Theta = \nabla F(I)$. We have four possible choices for the sign of the first order term of $F(I)$.
To fix ideas, we choose $  \frac{\pi}{2}\left(- I_1 - I_2\right)$ that is equivalent to $\frac{3\pi}{2}\left( I_1 + I_2\right)$.
$\Theta(I)$ has to satisfy the energy level condition for Highways in the reduced Poincar\'{e} function
\begin{align}\label{eq:original_eq_asymptotic}
\sum\limits_{k=1}^2A_k\cos(\Theta_k - \omega_k\tau^*(I,\Theta))\nonumber\\
+ A_3\left(\cos(-\tau^*(I,\Theta)) - 1\right) = 0,
\end{align}
and $\tau^*(I,\Theta)$ has to satisfy the equation of the crest
\begin{align}\label{eq:asymptotic_crests}
& \sum\limits_{k=1}^2\omega_kA_k\sin(\Theta_k - \omega_k\tau^*(I,\Theta))\nonumber\\
&+ A_3\sin(-\tau^*(I,\Theta))=0.
\end{align}

We want to write their version for $I_1$ and $I_2$ close to infinity.
Using \eqref{eq:lagrange_candidate} we notice that $\Theta_i = \Theta_i(I)$ takes the form
\begin{equation*}
\Theta_i = \frac{3\pi}{2} - a_i\sinh(\pi/2)\omega_i^3e^{\frac{-\pi\omega_i}{2}} + \mathcal{O}\left(\omega_1^2\omega_2^2e^{\frac{-\pi(\omega_1 + \omega_2)}{2}}\right).
\end{equation*}
This implies
\begin{align*}
\cos(\Theta_i - \omega_i\tau^*) &= -a_i\sinh(\pi/2)\omega_i^3e^{\frac{-\pi\omega_i}{2}}-\omega_i\tau^*_{\infty}\\
&+\mathcal{O}\left(\omega_1^2\omega_2^2e^{\frac{-\pi(\omega_1 + \omega_2)}{2}}\right)
\end{align*}
and
\begin{equation*}
\sin(\Theta_i - \omega_i\tau^*) = -1 + \mathcal{O}(\omega_i^6e^{-\pi\omega_i}),
\end{equation*}
In addition,
\begin{align*}
\cos(-\tau^*(I,\Theta)) &= 1-\frac{\tau^{*2}_{\infty}}{2} + \mathcal{O}(\tau^{*4}_{\infty}),\\
 \sin(-\tau^* ) &= -\tau^*_{\infty} + \mathcal{O}(\tau^*_{\infty}),
\end{align*}
where $\tau^*_{\infty}$ is an asymptotic approximation of $\tau^*$ that we will estimate in the following.
First, we notice that functions $A_1$ and $A_2$ can be approximated by
\begin{align*}
A_i&= 4\pi a_i\omega_i e^{\frac{-\pi\omega_i}{2}}\left( 1 + e^{-2\pi\omega_i} + \cdots\right)\\
& = 4\pi a_i\omega_i e^{\frac{-\pi\omega_i}{2}} + \mathcal{O}\left(\omega_ie^{\frac{-5\pi\omega_i}{2}}\right).
\end{align*}
From \eqref{eq:crests_cms}, the function $\tau^*(I,\Theta)$ satisfies
\begin{align*}
-\tau^*(I,\Theta) = -\arcsin\left(\sum\limits_{k=1}^2\frac{A_k\omega_k}{A_3}\sin(\Theta_k - \omega_k \tau^*(I,\Theta) ) \right),
\end{align*}
and therefore,
\begin{equation*}
\tau^*_{\infty} =\sum_{i=1}^2- 2a_i\sinh(\pi/2)\omega_i^2e^{\frac{-\pi\omega_i}{2}} + \mathcal{O}\left(\omega_ie^{\frac{-5\pi\omega_i}{2}}\right).
\end{equation*}

Applying these estimates in Eq.~\eqref{eq:original_eq_asymptotic} we obtain that the left-hand side of Eq.~\eqref{eq:original_eq_asymptotic} satisfies
\begin{align*}
&\sum_{i = 1}^2\left\{4\pi a_i \omega_ie^{-\pi\omega_1/2} \left[-a_i\sinh(\pi/2)\omega_i^3e^{-\pi\omega_i/2}\right.\right.\\
&\left.\left.- \omega_i\left(\sum_{k=1}^2 2a_k\sinh(\pi/2)\omega_k^2e^{-\pi\omega_k/2} \right)\right]\right\}\\
&-\frac{A_3}{2}\left(\sum_{i=1}^2 2a_i\sinh(\pi/2)\omega_i^2e^{-\pi\omega_i/2} \right)^2 \\
&+\mathcal{O}(\omega_1^2\omega_2^2e^{-\pi(\omega_1 + \omega_2)/2}) = \mathcal{O}(\omega_1^2\omega_2^2e^{-\pi(\omega_1 + \omega_2)/2}).
\end{align*}

In the same way, by applying in Eq.~\eqref{eq:asymptotic_crests} the estimates obtained, we have that the lef-hand side of Eq.~\eqref{eq:asymptotic_crests} satisfies
\begin{align*}
&-\sum_{i=1}^2 4\pi a_i\omega_i^2e^{\frac{-\pi \omega_i}{2}} +A_3\left(\sum_{i=1}^2 2a_i\sinh(\pi/2)\omega_i^2e^{\frac{-\pi\omega_i}{2}}\right)\\
& + \mathcal{O}(\omega_1^2\omega_2^2e^{-\pi(\omega_1 + \omega_2)/2})=\mathcal{O}(\omega_1^2\omega_2^2e^{-\pi(\omega_1 + \omega_2)/2}).
\end{align*}

Therefore, up to the order $\mathcal{O}(\omega_1^2\omega_2^2e^{-\pi(\omega_1 + \omega_2)/2})$, the equation of the crest and the energy level of the reduced Poincar\'{e} function are satisfied.
\end{proof}

So far, we could ensure the existence of Highways for some specific values of $I$.
However, in the exceptional case $a_1=a_2$, $\Omega_1=\Omega_2$, we can completely describe them and thus we can easily construct an explicit diffusing pseudo-orbit.

\begin{proposition}[Global explicit expression of Highways in a special case]
Consider the Hamiltonian \eqref{eq:ham}+\eqref{eq:pert} and $a_1 ~=~ a_2$ satisfying $2\left|a/a_3\right| < 0.625$ and $\Omega_1 = \Omega_2$.
Let $\phi_t^{-\mathcal{L}_0^*}(I^0,\Theta(I^0))$ be a scattering flow on a Highway such that
$I^0_1 = I^0_2$ and $\theta^0_1 = \theta^0_2$.
Then, $I_1(t) = I_2(t) = \bar{I}(t)$ and $\theta_1(t) = \theta_2(t) = \bar{\theta}(t)$ for all $t\in\mathbb{R}$ and can be described by
\begin{equation*}
\bar{\theta}_{\text{h}}(\bar{I}) = \left\{\begin{matrix}
\arccos\left(\frac{(1 - f(\bar{I}))A_3}{A}\right) + \bar{\omega}\arccos(f(\bar{I}));\\
\arccos\left(\frac{(1 - f(\bar{I}))A_3}{A}\right) - \bar{\omega}\arccos(f(\bar{I}));
\end{matrix}\right.
\end{equation*}
for $\bar{I}\leq 0$ and $\bar{I} > 0$, respectively, or
\begin{equation*}
\bar{\theta}_{\text{H}}(\bar{I}) = \left\{\begin{matrix}
-\arccos\left(\frac{(1 - f(\bar{I}))A_3}{A}\right) - \bar{\omega}\arccos(f(\bar{I}));\\
-\arccos\left(\frac{(1 - f(\bar{I}))A_3}{A}\right) + \bar{\omega}\arccos(f(\bar{I}));
\end{matrix}\right.,
\end{equation*}
for $\bar{I}\leq 0$ and $\bar{I} > 0$, respectively,
where $A:= 2A_1= 2A_2$ and
\begin{align*}
f(\bar{I}):= \left\{\begin{matrix}
1 - \frac{A}{2A^2_3},\quad \bar{\omega} =\pm 1\\
\frac{\bar{\omega}^2A_3 - \sqrt{A^2_3 + (\bar{\omega}^2-1)\bar{\omega}^2A^2}}{(\bar{\omega}^2-1)A_3}\quad \bar{\omega} \neq\pm 1.
\end{matrix}
\right.
\end{align*}

\end{proposition}

\begin{proof}
Consider the scattering flow $\phi_t^{-\mathcal{L}_0^*}(I^0,\theta^0)$ defined by the differential equations in \eqref{eq:edo_SM_k=0}.
Assuming $ \Omega_1 = \Omega_2 =: \Omega$ and $a_1 = a_2 =:a$ and taking the initial conditions satisfying $I(0)= I^0$ and $\theta(0) = \theta^0$ where $I^0_1 = I^0_2$ and $\theta^0_1 = \theta^0_2 $, the solution $(I(t),\theta(t))$ of \eqref{eq:edo_SM_k=0} satisfies $\theta_1(t) = \theta_2(t)=:\bar{\theta}(t)$ and $I_1(t) = I_2(t) =:\bar{I}(t)$ for all $t$.

If the flow $\phi_t^{-\mathcal{L}_0^*}(I^0,\theta^0)$ lies on a Highway, it has to satisfy two equations: the equation of the crests given in \eqref{eq:def_crests} and
\begin{equation*}
\mathcal{L}^*_0(I,\theta) = A_3.
\end{equation*}

These equations can be rewritten as
\begin{align}\label{eq:to_find_formula_Highway}
A\cos(\bar{\theta} - \bar{\omega}\tau^*(\bar{I},\bar{\theta})) + A_3\cos(-\tau^*(\bar{I},\bar{\theta})) = A_3\\
\bar{\omega}A\sin(\bar{\theta} - \bar{\omega}\tau^*(\bar{I},\bar{\theta})) + A_3 \sin(-\tau^*(\bar{I},\bar{\theta})) =0,\nonumber
\end{align}
where $\bar{\omega} := \omega_1 = \omega_2$ and $A:= A(\bar{\omega}) = 4\pi\bar{\omega}a/\sinh(\pi\bar{\omega}/2)$ for $\bar{\omega}\neq 0$ and $A = 8a$.

Multiplying by $\bar{\omega}$ the first equation in \eqref{eq:to_find_formula_Highway} we obtain
\begin{eqnarray*}
\bar{\omega}A\cos(\bar{\theta} - \bar{\omega}\tau^*(\bar{I},\bar{\theta})) = - \bar{\omega} A_3\left(\cos(-\tau^*(\bar{I},\bar{\theta}))-1\right) \\
\bar{\omega}A\sin(\bar{\theta} - \bar{\omega}\tau^*(\bar{I},\bar{\theta})) =- A_3\sin(-\tau^*(\bar{I},\bar{\theta})).
\end{eqnarray*}

We sum these two equations squared and we obtain
\begin{equation*}
\bar{\omega}^2A^2 = \left[\bar{\omega}A_3\left(\cos(-\tau^*(\bar{I},\bar{\theta}))-1\right)\right]^2 + A^2\sin^2(-\tau^*(\bar{I},\bar{\theta})).
\end{equation*}

After some arithmetical manipulations, we obtain the following equation of second degree in $\cos(-\tau^*(\bar{I},\bar{\theta}))$
\begin{align*}
(\bar{\omega}^2-1)A^2_3\cos^2(-\tau^*(\bar{I},\bar{\theta})) - 2\bar{\omega}^2A^2_3\cos(-\tau^*(\bar{I},\bar{\theta}))\\
+A^2_3(\bar{\omega}^2 + 1) - \bar{\omega}^2A^2 = 0.
\end{align*}

For $\bar{\omega} =\pm 1$, we have
\begin{align*}
\cos(-\tau^*(\bar{I},\bar{\theta})) = 1 - \frac{A}{2A^2_3}.
\end{align*}
Otherwise,
\begin{align*}
&\cos(-\tau^*(\bar{I},\bar{\theta})) =\\
& \frac{2\bar{\omega}^2A^2_3 \pm \sqrt{4\bar{\omega}^4A^4_3- 4(\bar{\omega}^2-1)A^2_3\left[A^2_3(\bar{\omega}^2 + 1) - \bar{\omega}^2A^2\right]}}{2(\bar{\omega}^2-1)A^2_3}.
\end{align*}
After more arithmetical manipulation and considering $-1\leq\cos(-\tau^*(\bar{I},\bar{\theta}))\leq 1$, we have
\begin{equation*}
\cos(-\tau^*(\bar{I},\bar{\theta})) = \frac{\bar{\omega}^2A_3 - \sqrt{A^2_3 + (\bar{\omega}^2-1)\bar{\omega}^2A^2}}{(\bar{\omega}^2-1)A_3}.
\end{equation*}
To simplify the notation, we define
\begin{align*}
f(\bar{I}):= \left\{\begin{matrix}
1 - \frac{A}{2A^2_3},\quad \bar{\omega} =\pm 1\\
\frac{\bar{\omega}^2A_3 - \sqrt{A^2_3 + (\bar{\omega}^2-1)\bar{\omega}^2A^2}}{(\bar{\omega}^2-1)A_3}\quad \bar{\omega} \neq\pm 1.
\end{matrix}
\right.
\end{align*}
And therefore,
\begin{equation*}
\Rightarrow  -\tau^*(\bar{I},\bar{\theta}) = \pm\arccos(f(\bar{I})).
\end{equation*}
We have two Highways.
This explains why we have found two different values for the function $\tau^*$.
Then we can rewrite the first equation of \eqref{eq:to_find_formula_Highway} as follows.
\begin{equation*}
A\cos(\bar{\theta} \pm \bar{\omega}\arccos(f(\bar{I}))) + A_3f(\bar{I}) = A_3
\end{equation*}
This immediately implies
\begin{equation*}
\bar{\theta} = \pm \arccos\left( \frac{ A_3\left(1 - f(\bar{I})\right)}{A} \right) \mp \bar{\omega}\arccos(f(\bar{I})).
\end{equation*}

From the four possibilities, by comparing with numerical results, we obtain that the Highways are described by
\begin{equation*}
\bar{\theta}_{\text{h}}(\bar{I}) = \left\{\begin{matrix}
\arccos\left(\frac{A_3(1 - f(\bar{I}))}{A}\right) + \bar{\omega}\arccos(f(\bar{I}));\\
\arccos\left(\frac{A_3(1 - f(\bar{I}))}{A}\right) - \bar{\omega}\arccos(f(\bar{I}));
\end{matrix}\right.
\end{equation*}
for $\bar{I}\leq 0$ and $\bar{I} > 0$, respectively, and
\begin{equation*}
\bar{\theta}_{\text{H}}(\bar{I}) = \left\{\begin{matrix}
-\arccos\left(\frac{A_3(1 - f(\bar{I}))}{A}\right) - \bar{\omega}\arccos(f(\bar{I})),;\\
-\arccos\left(\frac{A_3(1 - f(\bar{I}))}{A}\right) + \bar{\omega}\arccos(f(\bar{I}));
\end{matrix}\right.
\end{equation*}
for $\bar{I}\leq 0$ and $\bar{I} > 0$, respectively.
\end{proof}

\begin{remark}
A very similar result was proved for 2$+$1/2 d.o.f.\cite{RGS18}.
\end{remark}

Note that since an orbit on a Highway has the form $(I,\Theta(I))$, we have
\begin{align*}
\dot{\Theta} = D\Theta(I)\dot{I}.
\end{align*}
Therefore, $\dot{I} = 0$ implies $\dot{\Theta}=0$.
Thus, an equilibrium point in the plane $(I_1,I_2)$ is related to an equilibrium point of the scattering flow.
From Prop.~\ref{lem:fixed_points}, we know that the scattering flow has only four equilibrium points.
However, none of them lie on a Highway since $\mathcal{L}^*_0(\mathbf{p}) \neq A_3$, $\mathbf{p}=(0,0,0,0),\,(0,0,\pi,\pi),\,(0,0,\pi,0)$ and $(0,0,0,\pi)$.
We can conclude that there are no equilibrium points for Highways on $(I_1,I_2)$.

\begin{theorem}[Global Existence of Highways]
\label{thm:global_existence_highways}
Assume $a_1a_2a_3 \neq 0$ and $\left|a_1/a_3\right| + \left|a_2/a_3\right| < 0.625$ in Hamiltonian \eqref{eq:ham}$+$\eqref{eq:pert}.
Inside a Highway there are no orbits that are contained in a compact set on the plane $(I_1,I_2)$.
In particular, given any $0<c<C$,
there exists an orbit of the scattering flow $\phi_t^{-\mathcal{L^*_0}}(I,\theta)$ on a Highway and a time $t^*$ such that
$$ \left| I(t_0) \right| <c \quad \text{ and }\quad \left| I(t^*)\right|>C.$$
\end{theorem}

\begin{proof}
Suppose by contradiction that there is an orbit $(\gamma,\Theta(\gamma))$ on a Highway that is contained in a compact set on the plane $(I_1,I_2)$.
This implies that its $\omega$-limit set $w(\gamma)$ on this plane is nonempty. Since there are no equilibrium points on $(I_1,I_2)$, $w(\gamma)$ is formed by regular points.
From the Poincar\'e-Bendixson theorem, $\gamma$ is a periodic orbit.
However, this immediately implies an equilibrium point in the interior of $\gamma$, a contradiction.
Observing that $(0,0,\pm \pi/2, \pm \pi/2)$ belongs to a Highway, the theorem is proved for any $c>0$.
\end{proof}

As an immediate consequence, we can ensure that the orbits are diffusing along a Highway.

\begin{corollary}
Assume $a_1a_2a_3 \neq 0$ and $\left|a_1/a_3\right| + \left|a_2/a_3\right| < 0.625$ in Hamiltonian \eqref{eq:ham}$+$\eqref{eq:pert}.
Given any $0<c_j<C_j$, $j = 1,2$,
there is at least an orbit $(I^i , \theta^i)_{0\leq i < N}$ of the scattering map $\mathcal{S}_0$ such that
$$ \left| I_j^0 \right| <c_j \quad \text{ and }\quad \left| I^N_j\right|>C_j, \qquad j = 1,2.$$
\end{corollary}

\begin{figure}[h]
\begin{center}
\includegraphics[angle=-90,width=0.4\textwidth]{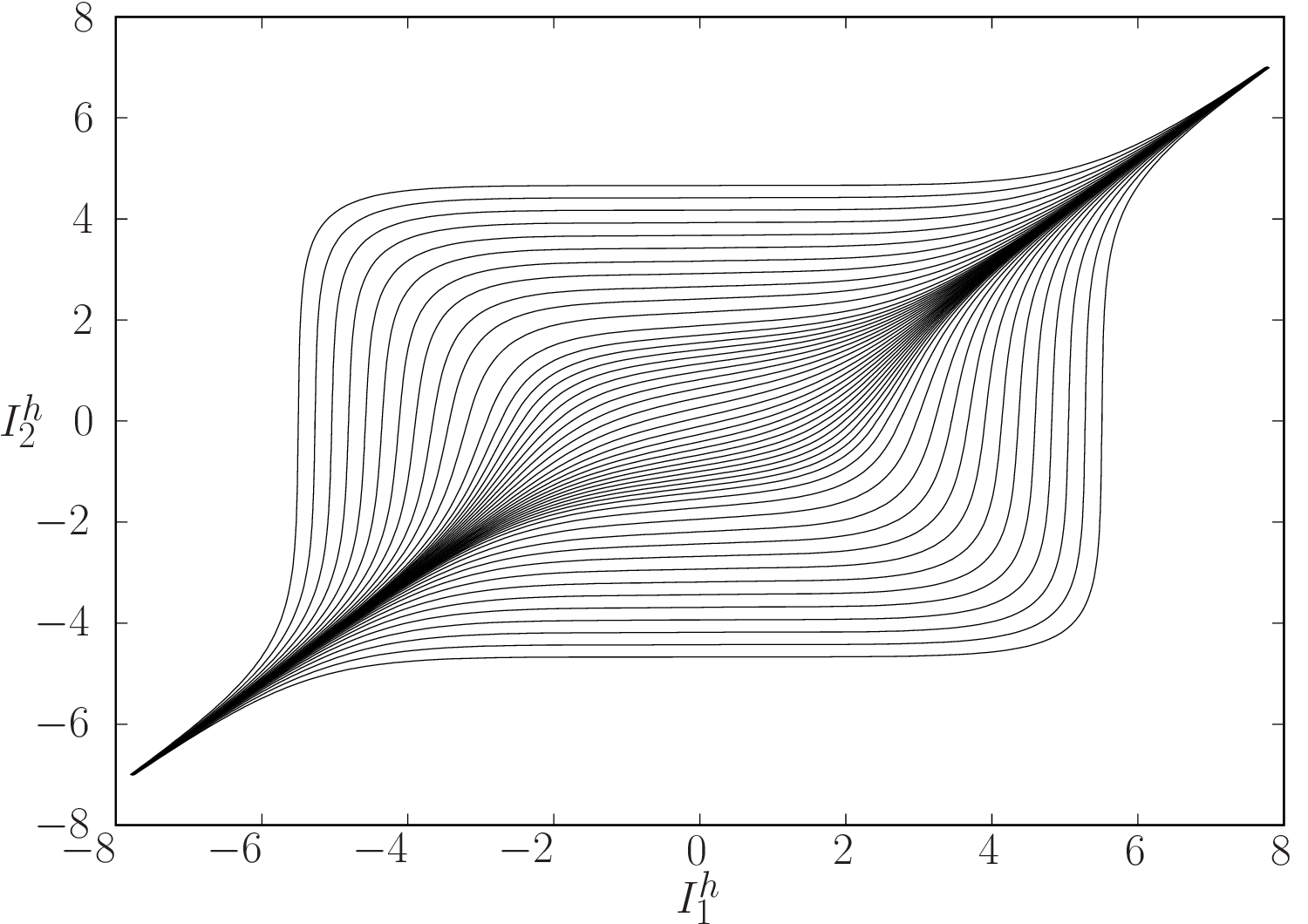}
\end{center}
\caption{Dynamics inside the Highway. Parameter values are $a_1=0.3$,
$a_2=0.1$, $a_3=1$ and $\Omega_1=\Omega_2=1$.}
\label{fig:dynamics_in_Highway}
\end{figure}

The following result shows that the orbits on a Highway do not wander along the plane $(I_1,I_2)$.
Indeed, they have an explicit asymptotic behavior, as we can see in Fig.~\ref{fig:dynamics_in_Highway}.

\begin{proposition}\label{prop:straight_line}
Assume $a_1a_2a_3 \neq 0$ and $\left|a_1/a_3\right| + \left|a_2/a_3\right| < 0.625$ in Hamiltonian \eqref{eq:ham}$+$\eqref{eq:pert}.
Let $(I^{\text{h}},\Theta(I^{\text{h}}))$ be a Highway.
For $I_2$, $I_1\gg 1$, we have
\begin{equation*}
I^{\text{h}}_2 = \frac{\Omega_1}{\Omega_2}I^{\text{h}}_1 - \frac{2}{\pi\Omega_2}\log\left(\frac{\Omega_1a_1}{\Omega_2a_2}\right) + \mathcal{O}(1/\omega_1),
\end{equation*}
and for $I_2$, $I_1\ll -1$,
\begin{equation*}
I^{\text{h}}_2 = \frac{\Omega_1}{\Omega_2}I^{\text{h}}_1 + \frac{2}{\pi\Omega_2}\log\left(\frac{\Omega_1a_1}{\Omega_2a_2}\right)+ \mathcal{O}(1/\omega_1),
\end{equation*}
\end{proposition}

\begin{proof}
Since $\dot{I}_i = -A_i\sin(\theta_i - \omega_i\tau^*)$ and $\omega_i = \Omega_i I_i$, we have
\begin{equation*}
\frac{d\omega_2}{d\omega_1} = \frac{\Omega_2A_2\sin(\theta_2 - \omega_2\tau^*)}{\Omega_1A_1\sin(\theta_1 - \omega_1\tau^*)}.
\end{equation*}
For $I_1, I_2\gg 1$, the above equation becomes
\begin{align*}
\frac{d\omega_2}{d\omega_1} & = \frac{4\pi a_2\Omega_2\omega_2 \exp(-\pi\omega_2/2) + \mathcal{O}(w_2\exp(-5\pi\omega_2/2))}{4\pi a_1\Omega_1\omega_e1 \exp(-\pi\omega_1/2)+ \mathcal{O}(w_1\exp(-5\pi\omega_1/2))}\\
& = \frac{a_2\Omega_2\omega_2\exp(-\pi\omega_2/2)+ \mathcal{O}(w_2\exp(-5\pi\omega_2/2))}{a_1\Omega_1\omega_1\exp(-\pi\omega_1/2)+ \mathcal{O}(w_1\exp(-5\pi\omega_1/2))},
\end{align*}
so that separating variables, we get
\begin{align*}
\frac{1}{a_2\Omega_2}\int\frac{\exp(\pi\omega_2/2)d\omega_2}{\omega_2} + \mathcal{O}\left(\frac{\exp\left(\frac{5\pi\omega_2}{2}\right)}{\omega_2^2}\right)= \\
\frac{1}{a_1\Omega_1}\int\frac{\exp(\pi\omega_1/2)d\omega_1}{\omega_1} + \mathcal{O}\left(\frac{\exp\left(\frac{5\pi\omega_1}{2}\right)}{\omega_1^2}\right),
\end{align*}
which can be also written as
\begin{align}
\frac{\exp(\pi\omega_2/2)}{a_2\Omega_2\pi\omega_2}+ \mathcal{O}(\omega_2^{-2}\exp(\pi\omega_2/2)) = \frac{\exp(\pi\omega_1/2)}{a_1\Omega_1\pi\omega_1}\nonumber\\
+ \mathcal{O}(\omega_1^{-2}\exp(\pi\omega_1/2)) + c_0.
\label{eq:dominant_terms}
\end{align}
Assume that $\omega_2$ takes the form $\omega_2 = \omega_1 + g(\omega_1)$.
Plugging this formula into Eq.~\eqref{eq:dominant_terms}, we obtain
\begin{align*}
\frac{\exp(\pi g(\omega_1)/2)}{a_2\Omega_2(\omega_1 + g(\omega_1))} +  \mathcal{O}\left(\frac{\exp\left(\pi g(\omega_1)/2\right)}{\left(\omega_1+g(\omega_1)\right)^2}\right) = \\
\frac{1}{a_1\Omega_1\omega_1} + \frac{c_0}{\exp(\pi\omega_1/2)} +\mathcal{O}(\omega_1^{-2}).
\end{align*}
As we assume that $\omega_1$ is large, the term with $\omega_1^{-2}$ dominates the term with $\exp(-\pi\omega_1)$.
Therefore,
\begin{equation}\label{eq:g_equation}
\exp(\pi g(\omega_1)/2) = \left(\frac{a_2\Omega_2}{a_1\Omega_1}\right)\left(1 + \frac{g(\omega_1)}{\omega_1}\right) + \mathcal{O}(\omega_1^{-1}).
\end{equation}
By solving Eq.~\eqref{eq:g_equation}, the function $g$ takes the following form.
\begin{equation*}
g(\omega_1) = \frac{-2}{\pi}\,W_{n}\left( \frac{-\pi a_1\Omega_1\omega_1}{2a_2\Omega_2} \exp(-\pi\omega_1/2) \right) - \omega_1 + \mathcal{O}(\omega_1^{-1}),
\end{equation*}
where $n\in\mathbb{Z}$ and $W_{n}$ is the Lambert $W$ function\citep{Corless1996}.
The Lambert $W$ function can be written both around $z= 0$ and $z\rightarrow \infty$ as\citep{Corless1996}
\begin{align}\label{eq:lambert_function}
W_{n}(z) &= \text{Log}(z) - \log\text{Log}(z) \nonumber\\
&+ \sum_{n= 0}^{\infty}\sum_{m= 1}^{\infty}c_{nm}(\log\text{Log}(z))^m(\text{Log}(z))^{-n-m},
\end{align}
where the coefficients $c_{nm}$ can be found using the Lagrange Inverse Theorem, $\text{Log}(z)$ is the complex logarithm for any no-principal branch, and $\log(z)$ is the complex logarithm for the principal branch.
As we are working with real numbers, Eq.~\eqref{eq:lambert_function} implies that the function $g$ is
\begin{align*}
&g(\omega_1) = \frac{-2}{\pi}\left(\log\left(\frac{\pi}{2}\left|\frac{a_1\Omega_1}{a_2\Omega_2} \right|\omega_1\exp\left(-\frac{\pi\omega_1}{2}\right)\right)\right.\\
&\left.- \log\left|\log\left(\frac{\pi}{2}\left|\frac{a_1\Omega_1}{a_2\Omega_2}\right|\omega_1 \exp\left(-\frac{\pi\omega_1}{2}\right)\right)\right| + \mathcal{O}\right)- \omega_1 + \mathcal{O},
\end{align*}
where $\mathcal{O}= \mathcal{O}(1/\omega_1)$.
After some algebraic manipulations, the function $g$ becomes
\begin{align*}
g(\omega_1) &= \frac{-2}{\pi}\left(\log\left(\frac{\pi}{2}\left|\frac{\pi a_1\Omega_1\omega_1}{a_2\Omega_2}\right| \right)\right.\\
&\left. - \log\left|\left(\log\left(\frac{\pi}{2}\left|\frac{a_1\Omega_1\omega_1}{a_2\Omega_2}\right| \right) -\pi\omega_1/2\right)\right|\right) +\mathcal{O}\left(\frac{1}{\omega}\right).
\end{align*}
We conclude that
\begin{equation*}
\lim_{\omega_1\rightarrow \infty} g(\omega_1) =\frac{2}{\pi}\log\left(\left| \frac{a_2\Omega_2}{a_1\Omega_1}\right|\right) + +\mathcal{O}\left(\frac{1}{\omega}\right).
\end{equation*}

For $\omega_2, \omega_1 \ll -1$ is analogous.

\end{proof}

\begin{figure}
\begin{center}
\includegraphics[angle=-90,width=0.4\textwidth]{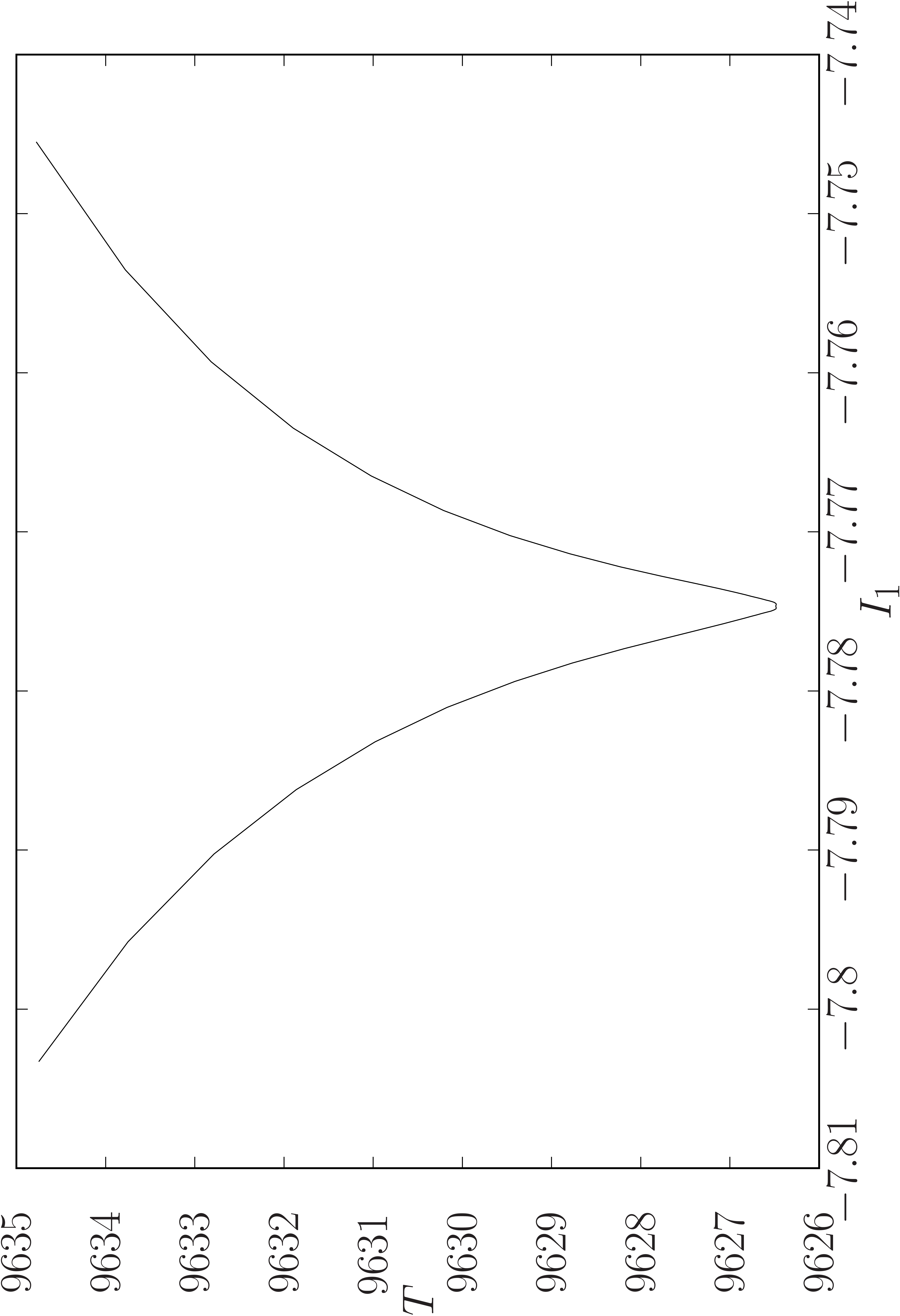}
\end{center}
\caption{Times taken by trajectories in
Fig.~\ref{fig:dynamics_in_Highway} to reach the section $I_2=7$ from
section $I_2=-7$.}
\label{fig:Highway_times}
\end{figure}

In Figure~\ref{fig:dynamics_in_Highway}, we simulate the dynamics of the scattering flow restricted to a Highway.
This simulation is carried out by taking initial
conditions along the sections given by $I_2=7$ and $I_2=-7$ (varying
$I_1$) and integrating backward and forward (respectively) so that the
trajectories match at $I_2=0$. The initial angles $\theta_1$ and
$\theta_2$ are taken by Equation~\eqref{eq:lagrange_candidate}.
Note that for large values of $I$, the dynamics inside the Highway is given by
the straight lines given in Prop.~\ref{prop:straight_line}.

In Figure~\ref{fig:Highway_times} we show the times taken by
these trajectories to travel from section $I_2=-7$ to section
$I_2=7$.
We can see that some orbits are faster than others; however, the difference in time is not quite significant.

\begin{remark}Fig.\ref{fig:dynamics_in_Highway} was made in \texttt{C++} by numerical integration of (28). The numerical method used is Runge-Kutta-Fehlberg (RKF78) with automatic step size control.
\end{remark}

\subsection{The time of diffusion}

We have two different types of time to estimate: the time of the scattering flow associated to the Hamiltonian $-\mathcal{L}_0^*(I,\theta)$ and the time along the homoclinic manifold to the NHIM $\tilde{\Lambda}$.
We will estimate this time only for orbits close to a Highway $(I,\Theta(I))$.
Recall that along the Highways, for large $I$ we can write $I_2$ asymptotically as a function of $I_1$, that is, there exists a function $f$ such that $I_2 = f(I_1)$.
Therefore, a diffusion from $I_1^0$ to $I_1^{\text{f}}$, both in $J$, immediately implies a diffusion from $I_2(I_1^0)$ to $I_2(I_1^{\text{f}})$.
This implies that we only need to study the time for the flow related to the differential equation for $I_1$.

The time of diffusion of an orbit close to a pseudoorbit built via iterates of scattering maps is estimated to be:
\begin{equation*}
T_{\text{d}} = N_{\text{s}}T_{\text{h}} + C_{\text{s}}T_{\text{i}},
\end{equation*}
where:	

\begin{itemize}
\item $N_{\text{s}}$ is the total number of iterates under scattering map. It is estimated by $T_{\text{s}}/\varepsilon$, where $T_{\text{s}}$ is the time that the scattering flow spends going through $\Delta I$. This comes from the fact that the scattering map jumps $\mathcal{O}(\varepsilon)$ along the level surfaces of $\mathcal{L}_0^*(I,\theta)$.

\item $T_{\text{h}}$ is the time under the flow along the homoclinic invariant manifolds of $\tilde{\Lambda}_{\varepsilon}$. This is the time spent by each application of the scattering map following the concrete homoclinic orbit to $\tilde{\Lambda}_{\varepsilon}$.

\item $C_{\text{s}}$ is the number of times that we need to apply the inner map to control the accumulated error, see Remark~\ref{rem:error}.

\item $T_{\text{i}}$ is the time under the inner map. This time appears if we use the inner map between iterates of the scattering map (it is sometimes called ergodization time).
\end{itemize}

We can refer to $N_{\text{s}}T_{\text{h}}$ as the time under scattering maps or simply \emph{outer time}, and to $ C_{\text{s}}T_{\text{i}}$ as the time under the inner map,
or simply \emph{inner time}.
In the following theorem, we will prove that the outer time is much larger than the inner time: $N_{\text{s}}T_{\text{h}} \gg C_{\text{s}}T_{\text{i}}$, for orbits close to Highways.
Hence, we prove that $T_{\text{d}} \approx N_{\text{s}}T_{\text{h}}$.

\begin{theorem}
%\label{prop:time}
The time of diffusion $T_{\text{d}}$ close to a Highway of Hamiltonian \eqref{eq:ham}$+$\eqref{eq:pert} with $\left|a_1/a_3\right| + \left|a_2/a_3\right| <0.625$ between $I_1^0$ and $I_1^{\text{f}}$ satisfies the following asymptotic expression.
\begin{equation}\label{eq:formula_time}
T_{\text{d}} = \frac{T_{\text{s}}}{\varepsilon}\left[2\log\left(\frac{C}{\varepsilon}\right) + \mathcal{O}(\varepsilon^{b})\right],
\end{equation}
for $\varepsilon\rightarrow 0$, where $0<b<1$, with
\begin{equation}\label{eq:t_s_diffusing_time}
T_{\text{s}} = \frac{1}{2\pi a_1\Omega_1}\int_{\omega_0}^{\omega_{\text{f}}}\frac{-\sinh(\pi\omega_1/2)d\omega_1}{\omega_1\sin(\theta_1 - \omega_1\tau^*)},
\end{equation}
where $\omega_0 = \Omega_1I_1^0 $ and $\omega_{\text{f}} = \Omega_1 I_{\text{f}}$, and
\begin{align*}
C = 16\left(\left|a_1\right| + \sum_{k=1}^2\frac{2\left|a_3\mu_k\right|\sinh(\pi/2)\left|\mu_1\right|}{\pi\left[1 - 1.466(\left|\mu_1\right| + \left|\mu_2\right|)\right]}M(\omega_k) \right).
\end{align*}
where $\alpha(\omega_i)$ was given in \eqref{eq:def_alpha}, $\mu_i = a_{i}/a_{3}$ and $M(\omega_i) = \max\limits_{I_i\in \left[I_i^0,I_i^\text{f}\right]}\left|\omega_i - \alpha(\omega_i)\right|. $
\end{theorem}

\begin{proof}

We estimate the time of diffusion to be
 $$T_{\text{d}}\approx T_{\text{s}} T_{\text{h}}/\varepsilon,$$
where $T_{\text{s}}$ is the time under the scattering flow on the Highway and $T_{\text{h}}$ is the time along the homoclinic manifolds to $\tilde{\Lambda}$.
We begin by studying $T_{\text{s}}$.

The differential equation for $I_1$ is given by
\begin{equation*}
\frac{dI_1}{dt} = -A_1\sin(\theta_1 - \omega_1\tau^*),
\end{equation*}
where $\theta_1 = \Theta_1(I_1)$ and $\tau^* = \tau^*(I_1)$.
Then, we have
\begin{equation*}
T_{\text{s}} = \frac{1}{2\pi a_1\Omega_1}\int_{\omega_0}^{\omega_{\text{f}}}\frac{-\sinh(\pi\omega_1/2)d\omega_1}{\omega_1\sin(\theta_1 - \omega_1\tau^*)},
\end{equation*}
where $\omega_0 = \Omega_1I_1^0$ and $\omega_{\text{f}} = \Omega_1 I_1^{\text{f}}$.

A point $(p_{\varepsilon}(\tau),q_{\varepsilon}(\tau))\in\overline{\text{B}_{\delta}(0)}\cap W_{\varepsilon}^{s,u}(0)$ lies on the perturbed separatrices on the plane $(p,q)$ and is given by $(p_{\varepsilon}(\tau),q_{\varepsilon}(\tau))= (p_0(\tau), q_0(\tau))  + \mathcal{O}(\varepsilon)$, where
$(p_0(\tau), q_0(\tau)) = \left(2/\cosh \tau,4\arctan e^{ \tau}\right)$.
The point $(p_0(\tau), q_0(\tau))$ can be asymptotically approximated by
\begin{align*}
p_{0}(\tau)  = \frac{4}{e^{\left|\tau\right|}}\left(1 + \mathcal{O}(e^{-2\left|\tau\right|})\right)\quad\text{ and}\\ q_{0}(\tau) = \mp\frac{4}{e^{\left|\tau\right|}}\left(1 +\mathcal{O}(e^{-2\left|\tau\right|})\right)\mod 2\pi.
\end{align*}

Taking into account the starting and ending points on $\partial \text{B}_{\delta}(0,0)$ and $\tau_{\text{f}} = -\tau_{\text{i}}$, we have
\begin{equation*}
\frac{4\sqrt{2}}{e^{u}}\left(1 + \mathcal{O}(e^{-2u})\right) = \delta.
\end{equation*}
where $u = \left|\tau_{\text{i}}\right|,\, \left|\tau_{\text{f}}\right|$.
Therefore,
\begin{equation*}
u = \log\left[\frac{4\sqrt{2}}{\delta}\left(1 + \mathcal{O}(\delta^2)\right)\right] = \log\left(\frac{4\sqrt{2}}{\delta}\right) + \mathcal{O}(\delta^2).
\end{equation*}

Then, the time along the homoclinic manifolds can be given by
\begin{equation}\label{eq:time_homoclinic}
T_{\text{h}}= 2\log\left(\frac{4\sqrt{2}}{\delta}\right) + \mathcal{O}(\delta^2) + \mathcal{O}(\varepsilon),
\end{equation}
where $\delta$ is the distance between the \NH and the piece of the invariant manifold that we are using to calculate the time.
Then, we have to estimate the value of $\delta$, since it must be small enough to preserve the scattering map.

Recall that from $p^2/2 + \cos q -1 =0$ the Melnikov potential given in \eqref{eq:mel_pot_gen}, can be written as
\begin{align*}
\mathcal{L}(I,\varphi,s)&=\int_{-\infty}^{+\infty}\frac{p^2(\sigma)}{2}\left(\sum_{k=1}^2a_k\cos(\varphi_k + \omega_k\sigma) \right.\\
&\left.+ a_3\cos(s + \sigma) \right)d\sigma.
\end{align*}

This implies that the reduced Poincar\'e function \eqref{eq:reduced_poincare_function} is
\begin{align*}
\mathcal{L}^*_0(I,\theta)=\frac{1}{2}\int_{-\infty}^{+\infty}p^2(\sigma)\left(a_1\cos(\theta_1 - \omega_1\tau^* + \omega_1\sigma)\right.\\
\left. + a_2\cos(\theta_2 - \omega_2\tau^* + \omega_2\sigma)+ a_3\cos(-\tau^*+ \sigma) \right)d\sigma.
\end{align*}

As we are considering just a piece of the homoclinic invariant manifold $\delta$-close to $\tilde{\Lambda}$, we are going to approximate the above integration by integrating for a finite interval of time $\left[t_0,t_{\text{f}}\right]$ in such a way that we have
\begin{equation}\label{eq:eq_condition_time}
\left|\frac{\partial \mathcal{L}^*_0}{\partial \theta_1} - \left(\frac{\partial \mathcal{L}^*_0}{\partial \theta_1}\right)_{\delta}\right| < \varepsilon,
\end{equation}
where
\begin{align*}
\left(\frac{\partial \mathcal{L}^*_0}{\partial \theta_1}\right)_{\delta} = \frac{1}{2}\int_{t_0}^{t_{\text{f}}}\frac{\partial}{\partial \theta_1}p^2(\sigma)\left(a_1\cos(\theta_1 - \omega_1\tau^* + \omega_1\sigma)\right.\\
\left. + a_2\cos(\theta_2 - \omega_2\tau^* + \omega_2\sigma) + a_3\cos(-\tau^*+ \sigma) \right)d\sigma.
\end{align*}
To simplify the calculations, we assume $t_0 = - t_{\text{f}}$.
Then, we have
\begin{align*}
&\left|\frac{\partial \mathcal{L}^*_0}{\partial \theta_1} - \left(\frac{\partial \mathcal{L}^*_0}{\partial \theta_1}\right)_{\delta}\right| =\\
&\left| \int_{t_{\text{f}}}^{+\infty}p^2(\sigma)\left[-a_1\left(1-\omega_1\frac{\partial \tau^*}{\partial \theta_1}\right)\sin(\theta_1 - \omega_1\tau^* + \omega_1\sigma) +\right. \right.\\
& \left.\left. a_2\omega_2\frac{\partial \tau^*}{\partial \theta_1}\sin(\theta_2 - \omega_2\tau^* + \omega_2\sigma) + a_3\frac{\partial \tau^*}{\partial \theta_1}\sin(-\tau^*+ \sigma) \right]d\sigma\right|\\
\end{align*}

Using the equation of the crest given in \eqref{eq:crests_cms} and the Triangular inequality
\begin{align}\label{eq:dif_der_time}
&\left|\frac{\partial \mathcal{L}^*_0}{\partial \theta_1} - \left(\frac{\partial \mathcal{L}^*_0}{\partial \theta_1}\right)_{\delta}\right| \leq  \left|a_1\right| \int_{t_{\text{f}}}^{+\infty}p^2(\sigma)d\sigma\nonumber\\
&+\sum_{k=1}^2\left|a_3\mu_k\right|\int_{t_{\text{f}}}^{+\infty}p^2(\sigma)\left|\frac{\partial \tau^*}{\partial \theta_1}\right|\left|\omega_k - \alpha(\omega_k)\right|d\sigma \nonumber\\
&\leq\int_{t_{\text{f}}}^{+\infty}p^2(\sigma)d\sigma\left(\left|a_1\right|\right.\nonumber\\
&\left.+\sum_{k=1}^2\left|a_3\mu_k\right| \left|\frac{\partial \tau^*}{\partial \theta_1}\right|\max\limits_{I_k\in \left[I_k^0,I_k^\text{f}\right]}\left|\omega_k - \alpha(\omega_k)\right|\right)d\sigma.
\end{align}

As $(I,\Theta(I))$ is on a Highway and from the equation of the crest, we have
\begin{equation*}
\frac{\partial\tau^*}{\partial \theta_1} = \frac{\alpha(\omega_1)\mu_1\cos(\theta_1 - \omega_1\tau^*)}{1 +\sum\limits_{k=1}^2 \left(\omega_k^2 -1\right)\frac{A_k}{A_3}\cos(\theta_k - \omega_k\tau^*)}
\end{equation*}

Applying the fact that we are assuming $\left|\mu_1\right| + \left|\mu_2\right| <0.625$,

\begin{equation*}
\left|\frac{\partial\tau^*}{\partial \theta_1}\right| \leq \frac{2\sinh(\pi/2)\left|\mu_1\right|}{\pi\left[1 - 1.466(\left|\mu_1\right| + \left|\mu_2\right|)\right]}.
\end{equation*}
In addition, we have $p(\sigma) = 4e^{-\left|\sigma\right|}\left(1 + \mathcal{O}(e^{-2\left|\sigma\right|})\right)$.
Therefore, inequality \eqref{eq:dif_der_time} can be rewritten as
\begin{align*}
\left|\frac{\partial \mathcal{L}^*_0}{\partial \theta_1} - \left(\frac{\partial \mathcal{L}^*_0}{\partial \theta_1}\right)_{\delta}\right|&\leq Ce^{-t_{\text{f}}} + \mathcal{O}(e^{-3t_{\text{f}}}),
\end{align*}
where
\begin{align}\label{eq:constant_time}
C = 16\left(\left|a_1\right| + \sum_{k=1}^2\frac{2\left|a_3\mu_k\right|\sinh(\pi/2)\left|\mu_1\right|}{\pi\left[1 - 1.466(\left|\mu_1\right| + \left|\mu_2\right|)\right]}M(\omega_k) \right),
\end{align}
with $M(\omega_i) = \max\limits_{I_i\in \left[I_i^0,I_i^\text{f}\right]}\left|\omega_i - \alpha(\omega_i)\right|. $
So, by using the expression of $T_{\text{h}}$ given in \eqref{eq:time_homoclinic}, we have
\begin{align*}
\left|\frac{\partial \mathcal{L}^*_0}{\partial \theta_1} - \left(\frac{\partial \mathcal{L}^*_0}{\partial \theta_1}\right)_{\delta}\right|&\leq \frac{C\delta (1 + \mathcal{O}(\delta^2))}{4\sqrt{2}},
\end{align*}

To satisfy Eq.~\eqref{eq:eq_condition_time}, we have to take a $\delta$ satisfying
\begin{equation*}
\delta = \frac{4\varepsilon\sqrt{2}}{C}(1 + \mathcal{O}(\varepsilon^2)).
\end{equation*}

Therefore, by inserting this value of $\delta$, time $T_{\text{h}}$ can be expressed as
\begin{equation*}
T_{\text{h}} = 2\log\left(\frac{C}{\varepsilon}\right) + \mathcal{O}(\varepsilon),
\end{equation*}
where $C$ is given by \eqref{eq:constant_time}.

Note that this time is only related to the outer time under the scattering map.
However, the inner map plays a role in the diffusion dynamics, and in principle the inner time should be taken into account.
As in the case of 2$+$1/2 d.o.f, this is not true\cite{Delshams2017}, since the inner time is much shorter than the outer, so it can be neglected.

We want to estimate the diffusion time for an orbit close to a Highway.
Then, let $(I^{\text{h}},\theta^{\text{h}})$ be a point on a Highway, take a point $(I^{\text{h}},\theta^{\text{h}}) + \Delta (I,\theta)$, with $\left\|\Delta (I,\theta)\right\|$ small enough.
From Gr\"{o}nwall's inequality, we have
\begin{align*}
\left\| \phi_t^{-\mathcal{L}^*_0}((I^{\text{h}},\theta^{\text{h}}) + \Delta (I,\theta)) -  \phi_t^{-\mathcal{L}^*_0}((I^{\text{h}},\theta^{\text{h}}) ) \right\| \leq\\  \left\|\Delta (I,\theta)) \right\|e^{K \left|t- t_0\right|}
\end{align*}
where $ K= \max_{s \in \left[t_0 , t\right]} \left\|\text{Hess}(\phi_s^{-\mathcal{L}^*_0}((I^{\text{h}},\theta^{\text{h}}))\right\|$ and $\text{Hess}$ is the Hessian matrix of $-\mathcal{L}^*_0$, and $\phi_t^{-\mathcal{L}^*_0}$ is the scattering flow.
To ensure that this orbit is close to the Highway we will take $\left\|\Delta(I,\theta)\right\|$ and $\left|t - t_0\right|$ small,
Then, $\left\|\Delta(I,\theta)\right\| = \varepsilon^{a}$ and $\left|t - t_0\right| =\varepsilon^c$ where $c$ is taken in such a way $e^{K \left|t- t_0\right|} = \mathcal{O}(1)$.
After this interval of time, we apply the inner map to reach a point closer to the Highway.
The time by this application of the inner map is estimated by applying the following theorem.
\begin{theorem}[Minkowski's theorem \cite{Siegel89}]\label{thm:mink}
For $Q>0$ and the linear operator $L:\mathbb{R}^{m}\rightarrow \mathbb{R}^{n}$ there exists a solution $(0,0)\neq (x,y) \in \mathbb{Z}^m\times\mathbb{Z}^n$ for the inequalities $\left\|x\right\|_{\infty}\leq Q$, $\left\|Lx - y\right\|_{\infty} \leq Q^{-\frac{m}{n}}$.
\end{theorem}

We consider that the inner dynamics is constant on the variable $I = (I_1, I_2)$, then we have to estimate $t$ such that $\left\|\varphi +t\omega - \varphi \mod 2\pi\right\| \leq \varepsilon^a$.
By considering a time 2$\pi$- periodic, this equation can be written as
\begin{equation*}
2\pi\left\| l\omega - k\right\| \leq \varepsilon^a,\quad l\in \mathbb{N} \text{ and } k\in \mathbb{Z}^{2}.
\end{equation*}
From Theorem~\ref{thm:mink}, the time $t$ satisfies $\left|t\right| \leq 2\pi\varepsilon^{-2a}$.
Therefore, we see that the time related to the inner dynamics $T_i$ satisfies
\begin{equation*}
\left|T_i\right| \leq \frac{T_s}{\varepsilon^{1+c}}\left(2\pi\varepsilon^{-2a}\right).
\end{equation*}

Observe that this time $T_i$ is comparable to the time under the scattering map if $\log\varepsilon^{-1}$ is comparable to $\varepsilon^{-c-2a}$.
But, taking $1 \gg -c \geq 2a >0$ we obtain $\log\varepsilon^{-1}\gg \varepsilon^{-c-2a}$.

We can now conclude that the diffusion time $T_d$ takes the form \eqref{eq:formula_time}, where $b = -c-2a$.
\end{proof}

\begin{remark}
We emphasize that the estimation of time obtained in the above theorem is very similar to the estimation for the diffusing time for 2$+$1/2 d.of. \cite{Delshams2017} and agrees with the estimation of `optimal' diffusing time\cite{Berti2003,Cresson2003,Tre04}. The main novelty here is that the constants $T_\text{s}$ and $C$ are explicitly given for diffusion orbits close to a Highway.
Note that $T_{\text{s}}$ depends on $\tau^*$, which does not have an analytical expression.
However, it can be easily computed using numerical methods.
\end{remark}
%\vspace*{0.1cm}
%\
%\appendix
%

% If in two-column mode, this environment will change to single-column format so that long equations can be displayed.
%\begin{widetext}
%$$\mbox{put long equation here}$$
%\end{widetext}

% Figures should be put into the text as floats.
% Use the graphics or graphicx packages (distributed with LaTeX2e).
% Put the label that you will use with \ref{} command in the braces of the \label{} command.
%
% \begin{figure}
% \end{figure}

% \begin{table}
% \end{tabular}
% \end{table}

\begin{acknowledgments}
AD has been partially supported by the Spanish MINECO/FEDER Grant PID2021-123968NB-I00.
RGS has been partially supported by CNPq, Conselho Nacional de Desenvolvimento Cient\'{i}fico e Tecnol\'{o}gico - Brasil and  the Priority Research Area SciMat under the program Excellence Initiative - Research University at the Jagiellonian University in Krak\'ow.
RGS would like to thank the hospitality of the Departament de Matem\`atiques of Universitat Polit\`ecnica de Catalunya and Matematiska institutionen of Uppsala Universitet where part of this work was carried out.
The authors would like to express their gratitude to the anonymous referees for their comments
and suggestions which have contributed to improve the final form of this paper.
\end{acknowledgments}

\bibliography{references}

\end{document}